\documentclass[a4paper,11pt]{amsart}
\usepackage{amssymb,amsmath}
\usepackage{amsopn}
\usepackage{mathrsfs}
\usepackage{color}
\usepackage{mathtools}
\usepackage{wasysym}
\usepackage{vmargin}
\usepackage{pdfsync}
\usepackage{vmargin}
\usepackage{color}
\usepackage{amscd}
\usepackage{verbatim}
\usepackage{bbm}
\newtheorem{theorem}{Theorem}[section]
\newtheorem{definition}[theorem]{Definition}
\newtheorem{lemma}[theorem]{Lemma}
\newtheorem{proposition}[theorem]{Proposition}
\newtheorem{corollary}[theorem]{Corollary}

\numberwithin{equation}{section}

\newcommand{\CC}{\mathbb{C}}

\newcommand{\rr}{\mathbb{R}}
\newcommand{\eps}{\varepsilon}
\newcommand{\nn}{\mathbb{N}}
\newcommand{\cc}{\mathbb{C}}
\def\un{{\mathrm{1~\hspace{-1.4ex}l}}}

\def\N{\mathbb N}
\def\Z{\mathbb Z}

\def\R{\mathbb R}

\def\val#1{\vert#1\vert}

\def\l2{L^2(\R^{n})}
\def\L2{L^2(\R^{2n})}
\def\supp{\operatorname{supp}}

\def\eps{\varepsilon}

\def\mat22#1#2#3#4{\begin{pmatrix}#1&#2\\ #3&#4\end{pmatrix}}

\author{J\'er\'emy \textsc{Martin} \&  Karel \textsc{Pravda-Starov}}

\address{\noindent \textsc{J\'er\'emy Martin, Univ Rennes, CNRS, IRMAR - UMR 6625, F-35000 Rennes, France
}}
\email{jeremy.martin@ens-rennes.fr}
\address{\noindent \textsc{Karel Pravda-Starov, Univ Rennes, CNRS, IRMAR - UMR 6625, F-35000 Rennes, France
}}
\email{karel.pravda-starov@univ-rennes1.fr}


\thanks{\textbf{Acknowledgements.} The authors express their gratefulness to the Centre de Math\'ematiques Henri Lebesgue for the very stimulating scientific environment.}

\begin{document}
\title[Spectral inequalities for Hermite functions and null-controllability]{Spectral inequalities for combinations of Hermite functions and null-controllability for evolution equations enjoying Gelfand-Shilov smoothing effects}
\begin{abstract}
This work is devoted to the study of uncertainty principles for finite combinations of Hermite functions. We establish some spectral inequalities for
control subsets that are thick with respect to some unbounded densities growing almost linearly at infinity, and provide quantitative estimates, with respect to the energy level of the Hermite functions seen as eigenfunctions of the harmonic oscillator, for the constants appearing in these spectral estimates. These spectral inequalities allow to derive the null-controllability in any positive time for evolution equations enjoying specific regularizing effects. More precisely, for a given index $\frac{1}{2} \leq \mu <1$, we deduce sufficient geometric conditions on control subsets to ensure the null-controllability of evolution equations enjoying regularizing effects in the symmetric Gelfand-Shilov space $S^{\mu}_{\mu}(\rr^n)$. These results apply in particular to derive the null-controllability in any positive time for evolution equations associated to certain classes of hypoelliptic non-selfadjoint quadratic operators, or to fractional harmonic oscillators.
\end{abstract}
\maketitle
\section{Introduction}

The classical uncertainty principle was established by Heisenberg and is linked to the impossibility to determine precisely the position and the momentum of quantum particles. Uncertainty principles are mathematical results that give limitations on the simultaneous concentration of a function and its Fourier transform. There are various uncertainty principles with formulations of different nature. A formulation of uncertainty principles is for instance that a non-zero function and its Fourier transform cannot both have small supports. In particular, a non-zero $L^2(\rr)$-function whose Fourier transform is compactly supported extends as a non-zero entire function with full support thanks to the isolated zeros theorem. Another formulation of uncertainty principles can be illustrated by the following notions of weak and strong annihilating pairs:

\medskip

\begin{definition} [Annihilating pairs]
Let $S,\Sigma$ be two measurable subsets of $\rr^n$. 
\begin{itemize}
\item[-] The pair $(S,\Sigma)$ is said to be a \emph{weak annihilating pair} if the only function $f\in L^2(\rr^n)$
with $\supp f\subset S$ and $\supp\widehat{f} \subset \Sigma$ is zero $f=0$
\item[-] The pair $(S,\Sigma)$ is said to be a \emph{strong annihilating pair} if there exists a positive constant $C=C(S,\Sigma)>0$ such that 
for all $f \in L^2(\rr^n)$,
\begin{equation}\label{strongly}
\int_{\rr^n}|f(x)|^2dx \leq C\Big(\int_{\rr^n \setminus S}|f(x)|^2dx  
+\int_{\rr^n \setminus \Sigma}|\widehat{f}(\xi)|^2d\xi\Big)
\end{equation}
\end{itemize}
\end{definition}

\medskip

It can be readily checked that a pair $(S, \Sigma)$ is a strong annihilating one if and only if there exists a positive constant $D=D(S, \Sigma)>0$ such that for all $f \in L^2(\rr^n)$ with $\supp\widehat{f} \subset \Sigma$, 
\begin{equation}\label{strongly2}
\|f\|_{L^2(\rr^n)} \leq D\|f\|_{L^2(\rr^n \setminus S)}.
\end{equation}
As already mentioned above, the pair $(S,\Sigma)$ is a weak annihilating one if $S$ and $\Sigma$ are compact sets.  More generally, Benedicks has shown in \cite{benedicks} that $(S,\Sigma)$ is a weak annihilating pair if $S$ and $\Sigma$ are sets of finite Lebesgue measure 
$|S|, |\Sigma| <+\infty$. Under this assumption, the result of Amrein-Berthier \cite{amrein} actually shows that the pair 
$(S,\Sigma)$ is a strong annihilating one. The estimate $C(S,\Sigma)\leq \kappa e^{\kappa |S||\Sigma|}$ (which is sharp up to the numerical 
constant $\kappa>0$) has been established by Nazarov~\cite{Na} in dimension $n=1$. This result was extended in the multi-dimensional case by Jaming~\cite{Ja}, with the quantitative estimate 
$$C(S,\Sigma)\leq \kappa e^{\kappa (|S||\Sigma|)^{1/n}},$$ 
holding if in addition one of the two subsets of finite Lebesgue measure $S$ or~$\Sigma$ is convex.

An exhaustive description of all strong annihilating pairs seems for now totally out of reach. We refer the reader for instance to the works~\cite{AO,BD1,BD2,De,DJ,SVW} for a large variety of results and techniques available, as well as for examples of weak annihilating pairs that are not strong annihilating ones.
On the other hand, there is an exhaustive description of all the support sets $S$ forming a strong annihilating pair with any bounded spectral set $\Sigma$. This description is given by the Logvinenko-Sereda theorem \cite{Logvinenko_Sereda}:

\medskip

\begin{theorem}[Logvinenko-Sereda]\label{sereda}
Let $S,\Sigma\subset\rr^n$ be measurable subsets with $\Sigma$ bounded.
The following assertions are equivalent:
\begin{itemize}
\item[-] The pair $(S,\Sigma)$ is a strong annihilating pair
\item[-] The subset $\rr^n\setminus S$ is thick, that is, there exist a cube $K \subset \rr^n$ with sides parallel to coordinate axes and a positive constant $0<\gamma \leq 1$ such that 
$$\forall x \in \rr^n, \quad |(K+x) \cap (\rr^n\setminus S)| \geq \gamma|K|>0,$$
where $|A|$ denotes the Lebesgue measure of the measurable set $A$.
\end{itemize}
\end{theorem}

\medskip

It is noticeable to observe that if $(S,\Sigma)$ is a strong annihilating pair for some bounded subset $\Sigma$, then $S$ makes up a strong annihilating pair with every bounded subset $\Sigma$, but the above constants $C(S,\Sigma)>0$ and $D(S,\Sigma)>0$ do depend on $\Sigma$.
In order to be able to use this result in the control theory of partial differential equations, it is then essential to understand how the positive constant $D(S,\Sigma)>0$ depends on the bounded set $\Sigma$. This question was addressed by Kovrijkine~\cite{Kovrijkine} (Theorem~3) who established the following quantitative estimates:

\medskip

\begin{theorem}[Kovrijkine] \label{Kovrij} There exists a universal positive constant $C_n>0$ depending only on the dimension $n\geq 1$ such that
if $\omega$ is a $\gamma$-thick set at scale $L>0$, that is, 
\begin{equation}\label{thick1v}
\forall x \in \rr^n, \quad |\omega \cap (x+[0,L]^n)| \geq \gamma L^n,
\end{equation}
with $0<\gamma \leq 1$,
then, for all $R>0$ and $f \in L^2(\rr^n)$ with 
$\supp\widehat{f} \subset [-R,R]^n$, the following estimate holds
\begin{equation}\label{kov}
\|f\|_{L^2(\rr^n)} \leq \Big(\frac{C_n}{\gamma}\Big)^{C_n(1+LR)}\|f\|_{L^2(\omega)}.
\end{equation}
\end{theorem}

\medskip

In all this work, the Fourier transform is used with the following normalization
$$\widehat{f}(\xi)=\int_{\rr^n}f(x)e^{-ix \cdot \xi}d\xi, \qquad \xi \in \rr^n.$$
Given a measurable subset, notice that it is thick in $\rr^n$ if and only if it is $\gamma$-thick at scale $L$ for some positive constants $0<\gamma \leq 1$ and $L>0$. Thus, the notion of $\gamma$-thickness at a positive scale allows to quantify the general thickness property.

Thanks to this explicit dependence of the constant with respect to the parameter $R>0$ in the estimate (\ref{kov}), Egidi  and Veseli\'c~\cite{veselic}, and Wang, Wang, Zhang and Zhang~\cite{Wang} have independently established that the heat equation 
\begin{equation}\label{heat}
\left\lbrace \begin{array}{ll}
(\partial_t -\Delta_x)f(t,x)=\un_{\omega}(x)u(t,x)\,, \quad &  x \in \mathbb{R}^n,\ t>0, \\
f|_{t=0}=f_0 \in L^2(\rr^n),                                       &  
\end{array}\right.
\end{equation}
is null-controllable in any positive time $T>0$ from a measurable control subset $\omega \subset \rr^n$ if and only if this subset $\omega$ is thick in $\rr^n$. The recent work~\cite{egidi} by Beauchard, Egidi and the second author has shown that this geometric necessary and sufficient condition on control subsets to ensure null-controllability, extends more generally for hypoelliptic non-autonomous Ornstein-Uhlenbeck equations when the moving control subsets comply with the flow associated to the transport part of the Ornstein-Uhlenbeck operators. 

The notion of null-controllability is defined as follows:

\medskip

\begin{definition} [Null-controllability] Let $P$ be a closed operator on $L^2(\rr^n)$ which is the infinitesimal generator of a strongly continuous semigroup $(e^{-tP})_{t \geq 0}$ on $L^2(\rr^n)$, $T>0$ and $\omega$ be a measurable subset of $\mathbb{R}^n$. 
The evolution equation 
\begin{equation}\label{syst_general}
\left\lbrace \begin{array}{ll}
(\partial_t + P)f(t,x)=\un_{\omega}(x)u(t,x), \quad &  x \in \mathbb{R}^n,\ t>0, \\
f|_{t=0}=f_0 \in L^2(\rr^n),                                       &  
\end{array}\right.
\end{equation}
is said to be {\em null-controllable from the set $\omega$ in time} $T>0$ if, for any initial datum $f_0 \in L^{2}(\mathbb{R}^n)$, there exists a control function $u \in L^2((0,T)\times\mathbb{R}^n)$ supported in $(0,T)\times\omega$, such that the mild \emph{(}or semigroup\emph{)} solution of \eqref{syst_general} satisfies $f(T,\cdot)=0$.
\end{definition}

\medskip

By the Hilbert Uniqueness Method, see \cite{coron_book} (Theorem~2.44) or \cite{JLL_book}, the null controllability of the evolution equation \eqref{syst_general} is equivalent to the observability of the adjoint system 
\begin{equation} \label{adj_general}
\left\lbrace \begin{array}{ll}
(\partial_t + P^*)g(t,x)=0, \quad & x \in \mathbb{R}^n, \ t>0, \\
g|_{t=0}=g_0 \in L^2(\rr^n),
\end{array}\right.
\end{equation}
where $P^*$ denotes the $L^2(\rr^n)$-adjoint of $P$. 
The notion of observability is defined as follows:

\medskip

\begin{definition} [Observability] Let $T>0$ and $\omega$ be a measurable subset of $\mathbb{R}^n$. 
The evolution equation \eqref{adj_general} is said to be {\em observable from the set $\omega$ in time} $T>0$, if there exists a positive constant $C_T>0$ such that,
for any initial datum $g_0 \in L^{2}(\mathbb{R}^n)$, the mild \emph{(}or semigroup\emph{)} solution of \eqref{adj_general} satisfies
\begin{equation}\label{eq:observability}
\int\limits_{\mathbb{R}^n} |g(T,x)|^{2} dx  \leq C_T \int\limits_{0}^{T} \Big(\int\limits_{\omega} |g(t,x)|^{2} dx\Big) dt\,.
\end{equation}
\end{definition}

\medskip

Following~\cite{veselic}, the necessity of the thickness condition for control subsets to ensure the null-controllability of the heat equation is a consequence of a quasimodes construction; whereas the sufficiency is derived from an abstract observability result based on an adapted Lebeau-Robbiano method established by Beauchard and the second author with some contributions of Miller in~\cite{KK1} (Theorem~2.1). This abstract observability result whose proof is inspired from the works~\cite{mi1,mi2},
was extended in~\cite{egidi} (Theorem~3.2) to the non-autonomous case with moving control supports and under weaker dissipation estimates allowing controlled blow-up for small times in the dissipation estimates. The following statement is a simplified formulation of Theorem~3.2 in~\cite{egidi} limited to the semigroup case with fixed control supports and weaker dissipation estimates than in~\cite{KK1} (Theorem~2.1):

\medskip

\begin{theorem}[Beauchard, Egidi \& Pravda-Starov] \label{Meta_thm_AdaptedLRmethod}
 Let $\Omega$ be an open subset of $\mathbb{R}^n$,
 $\omega$ be a measurable subset of $\Omega$,
 $(\pi_k)_{k \geq 1}$ be a family of orthogonal projections on $L^2(\Omega)$,
 $(e^{-tA})_{t \geq 0}$ be a strongly continuous contraction semigroup on $L^2(\Omega)$; 
 $c_1, c_2, c_1', c_2',a, b, t_0, m_1>0 $ be positive constants with $a<b$; $m_2 \geq 0$.
If the following spectral inequality
\begin{equation} \label{Meta_thm_IS}
\forall g \in L^2(\Omega), \forall k \geq 1, \quad \|\pi_k g \|_{L^2(\Omega)} \leq c_1' e^{c_1 k^a} \|\pi_k g \|_{L^2(\omega)},
\end{equation}
and the following dissipation estimate with controlled blow-up 
\begin{equation} \label{Meta_thm_dissip}
\forall g \in L^2(\Omega), \forall k \geq 1, \forall 0<t<t_0, \quad \| (1-\pi_k)(e^{-tA} g)\|_{L^2(\Omega)} \leq \frac{e^{-c_2 t^{m_1} k^b}}{c_2' t^{m_2}} \|g\|_{L^2(\Omega)},
\end{equation}
hold, then there exists a positive constant $C>1$ such that the following observability estimate holds
\begin{equation} \label{meta_thm_IO}
\forall T>0, \forall g \in L^2(\Omega), \quad \| e^{-TA} g \|_{L^2(\Omega)}^2 \leq C\exp\Big(\frac{C}{T^{\frac{am_1}{b-a}}}\Big) \int_0^T \|e^{-tA} g \|_{L^2(\omega)}^2 dt.
\end{equation}
\end{theorem}

\medskip

Notice that the assumptions in the above statement do not require that the orthogonal projections $(\pi_k)_{k \geq 1}$ are in any manner related to the spectral projections onto the eigenspaces of the infinitesimal generator $A$, which is allowed to be non-selfadjoint. According to the above statement, there are two key ingredients to derive a result of observability, or equivalently a result of null-controllability for the adjoint system, while using Theorem~\ref{Meta_thm_AdaptedLRmethod}, namely a spectral inequality (\ref{Meta_thm_IS}) and a dissipation estimate (\ref{Meta_thm_dissip}). For the heat equation, the orthogonal projections used are the frequency cutoff operators given by the orthogonal projections onto the closed vector subspaces
\begin{equation}\label{cutoff}
E_k=\big\{f \in L^2(\rr^n) : \textrm{ supp }\widehat{f} \subset [-k,k]^n\big\}, \qquad k \geq 1.
\end{equation}
With this choice, the dissipation estimate readily follows from the explicit formula
\begin{equation}\label{heat1}
\widehat{(e^{t\Delta_x}g)}(t,\xi)=\widehat{g}(\xi)e^{-t|\xi|^2}, \quad t \geq 0, \ \xi \in \rr^n,
\end{equation}     
whereas the spectral inequality is given by the sharpened formulation of the Logvinenko-Sereda theorem established by Kovrijkine (Theorem~\ref{Kovrij}). Notice that the power $1$ for the parameter $R$ in the estimate (\ref{kov}) and the power $2$ for the term $|\xi|$ in formula (\ref{heat1}) account for the fact that Theorem~\ref{Meta_thm_AdaptedLRmethod} can be applied with the parameters $a=1$, $b=2$ that satisfy the required condition $0<a<b$. It is therefore essential that the power of the parameter $R$ in the exponent of the estimate (\ref{kov}) is strictly less than $2$. Let us underline that Theorem~\ref{Meta_thm_AdaptedLRmethod} does not only apply with the use of frequency cutoff projections and a dissipation estimate induced by some Gevrey type regularizing effects. Other regularities than the Gevrey one can be taken into account. In this work, we are interested in obtaining results of null-controllability for evolution equations enjoying some regularizing effects in Gelfand-Shilov spaces. More specifically, given an abstract evolution equation enjoying some Gelfand-Shilov regularizing effects, we aim at finding sufficient geometric conditions on control subsets to ensure null-controllability in any positive time. 
The definition and basic properties related to Gelfand-Shilov regularity are recalled in appendix (Section~\ref{gelfand}). As recalled in this section, the Gelfand-Shilov regularity is characterized by specific exponential decays of both the functions and their Fourier transforms. In the symmetric case, the 
Gelfand-Shilov regularity can be read on the exponential decay of the Hermite coefficients when expanding the functions in the $L^2(\rr^n)$-Hermite basis $(\Phi_{\alpha})_{\alpha \in \nn^n}$. We refer the reader to Section~\ref{weighted} for the definition and some notations related to Hermite functions. Thanks to this second characterization of the Gelfand-Shilov regularity, a natural choice for the orthogonal projections $(\pi_k)_{k\geq 1}$ in order to apply Theorem~\ref{Meta_thm_AdaptedLRmethod} to prove the null-controllability of evolution equations enjoying some symmetric Gelfand-Shilov regularizing effects, are given by the Hermite orthogonal projections onto the closed vector subspaces in $L^2(\rr^n)$,
\begin{equation}\label{berk1}
\mathcal{E}_k=\textrm{Span}_{\cc}\{\Phi_\alpha\}_{\alpha \in \nn^n,\  |\alpha| \leq k}, \quad k \in \nn,
\end{equation}
where $\nn$ denotes the set of non-negative integers, and $|\alpha |= \alpha_1 +... + \alpha_n$, when $\alpha=(\alpha_1,...,\alpha_n) \in \nn^n$, that is, the orthogonal projections
\begin{equation}\label{projection}
\pi_k =\sum_{j=0}^k \mathbb{P}_j, \qquad \mathbb{P}_k g = \sum \limits_{\substack{\alpha \in \nn^n, \\ |\alpha|=k}} \left\langle g, \Phi_\alpha \right\rangle_{L^2(\rr^n)} \Phi_\alpha, \qquad k \geq 0,
\end{equation}
where $\mathbb{P}_k$ denotes the orthogonal projection onto the $k^{\text{th}}$ energy level associated with the harmonic oscillator
\begin{equation}\label{berk2}
\mathcal{H}= - \Delta_x + |x|^2= \sum_{k=0}^{+\infty} (2k+n) \mathbb{P}_k.
\end{equation}
Given an abstract evolution equation enjoying some symmetric Gelfand-Shilov regularizing effects, the dissipation estimate (\ref{Meta_thm_dissip}) is then expected to hold for the Hermite orthogonal projections $(\pi_k)_{k\geq 1}$ with some specific positive parameter $b>0$ related to the index of Gelfand-Shilov regularity. Let us notice that this dissipation estimate does not depend on the geometry of the control subset and that this geometry only plays a (key) role in the spectral inequality (\ref{Meta_thm_IS}). Adressing the problem of finding sufficient geometric conditions on control subsets to derive an observability result for this abstract evolution equation is therefore reduced to obtain quantitative spectral estimates of the type
\begin{equation}\label{eq7l}
\forall k \geq 1, \exists C_k(\omega)>0, \forall f \in L^2(\mathbb{R}^n), \quad \|\pi_k f\|_{L^2(\rr^n)} \leq C_k(\omega) \|\pi_k f\|_{L^2(\omega)},
\end{equation}
and figure out the largest class of control subsets for which the spectral inequality (\ref{Meta_thm_IS}) holds with some positive parameter $0<a<b$. This problem of studying 
under which conditions on the control subset~$\omega \subset \rr^n$, the spectral inequality (\ref{eq7l}) holds and how the geometric properties of the control subset $\omega$ relate to the possible growth of the positive constant $C_k(\omega)>0$ with respect to the energy level when $k \to +\infty$ was studied by Beauchard, Jaming and the second author in~\cite{kkj}. By a simple argument of equivalence of norms in finite dimension, the first result in~\cite{kkj} shows that for any measurable subset $\omega \subset \rr^n$ of positive Lebesgue measure $|\omega|>0$ and all $N \in \nn$, there does exist a positive constant $C_N(\omega)>0$ depending on $\omega$ and $N$ such that the following spectral inequality holds
\begin{equation}\label{spec}
\forall f \in \mathcal E_{N}, \quad \|f\|_{L^2(\rr^n)} \leq C_N(\omega)\|f\|_{L^2(\omega)}.
\end{equation} 
The main result in~\cite{kkj} (Theorem~2.1) then provides the following quantitative upper bounds on the positive constant $C_N(\omega)>0$ for the following three different geometries: 

\medskip

\noindent
$(i)$ If $\omega$ is a non-empty open subset of $\rr^n$, then there exists a positive constant $C=C(\omega)>1$ such that  
\begin{equation}\label{berk10}
\forall N \in \nn, \forall f \in \mathcal E_{N}, \quad \|f\|_{L^2(\rr^n)} \leq  Ce^{\frac{1}{2}N \ln(N+1)+CN}\|f\|_{L^2(\omega)}.
\end{equation}
$(ii)$ If the measurable subset $\omega \subset \rr^n$ satisfies the condition 
\begin{equation}\label{liminf}
\liminf_{R \to +\infty}\frac{|\omega \cap B(0,R)|}{|B(0,R)|}=\lim_{R \to +\infty}\Big(\inf_{r \geq R}\frac{|\omega \cap B(0,r)|}{|B(0,r)|}\Big)>0,
\end{equation}
where $B(0,R)$ denotes the open Euclidean ball in $\rr^n$ centered in $0$ with radius $R>0$, then there exists a positive constant $C=C(\omega)>1$ such that
\begin{equation}\label{berk11}
\forall N \in \nn, \forall f \in \mathcal E_{N}, \quad \|f\|_{L^2(\rr^n)} \leq  Ce^{CN}\|f\|_{L^2(\omega)}.
\end{equation}
$(iii)$ If the measurable subset $\omega \subset \rr^n$ is $\gamma$-thick at scale $L>0$, that is, when (\ref{thick1v}) holds, then there exist a positive constant $C=C(L,\gamma,n)>0$ depending on the dimension $n \geq 1$ and the parameters $0<\gamma \leq 1$, $L>0$, and a universal positive constant $\kappa=\kappa(n)>0$ only depending on the dimension such that 
\begin{equation}\label{berk12}
\forall N \in \nn, \forall f \in \mathcal E_{N}, \quad \|f\|_{L^2(\rr^n)} \leq  C\Big(\frac{\kappa}{\gamma}\Big)^{\kappa L\sqrt{N}}\|f\|_{L^2(\omega)}.
\end{equation}

These results show that the spectral inequality (\ref{Meta_thm_IS}) is satisfied with parameter $a=\frac{1}{2}$, when the control subset $\omega \subset \rr^n$ is $\gamma$-thick at scale $L>0$; whereas it holds with the parameter $a=1$ when the geometric condition (\ref{liminf}) holds.

The main result in the present work (Theorem~\ref{Spectral}) bridges the gap between the two spectral estimates (\ref{berk11}) and (\ref{berk12}) by figuring out sharp geometric conditions on the control subsets ensuring that the spectral inequality (\ref{Meta_thm_IS}) holds for any given parameter $\frac{1}{2} \leq a<1$. Given an abstract evolution equation enjoying some regularizing effects in the symmetric Gelfand-Shilov space $S_{\mu}^{\mu}(\rr^n)$, with $\frac{1}{2} \leq \mu<1$, some sharp sufficient geometric conditions on control subsets to ensure null-controllability are then deduced in Theorem~\ref{control_PDEfrac}, and some applications to derive the null-controllability of evolution equations associated to certain classes of hypoelliptic non-selfadjoint quadratic operators, or to fractional harmonic oscillators are given in Corollaries~\ref{th2_bis} and~\ref{corofrac}.

 \section{Statements of the main results}
 \subsection{Uncertainty principles for finite combinations of Hermite functions.}
The main result contained in the present work is the following uncertainty principles for finite combinations of Hermite functions:

\medskip

\begin{theorem}\label{Spectral}
Let $\rho : \rr^n \longrightarrow (0,+\infty)$ be a $\frac{1}{2}$-Lipschitz positive function with $\rr^n$ being equipped with the Euclidean norm,
satisfying that there exist some positive constants $0< \eps \leq 1$, $m>0$, $R>0$ such that
\begin{equation*}
\forall x \in \rr^n, \quad 0<m \leq \rho(x) \leq R{\left\langle x\right\rangle}^{1-\eps}.
\end{equation*}
Let $\omega$ be a measurable subset of $\rr^n$ which is $\gamma$-thick with respect to the density $\rho$, that is,
\begin{equation}{\label{thick_rho0}}
 \exists 0< \gamma \leq 1, \forall x \in \rr^n, \quad |\omega \cap B(x,\rho(x))| \geq \gamma \left|B(x,\rho(x))\right|,
\end{equation}
where $B(y,r)$ denotes the Euclidean ball centered at $y \in \rr^n$ with radius $r>0$, and $|\cdot|$ denotes the Lebesgue measure. 
Then, there exist some positive constant $\kappa_n(m, R, \gamma, \eps)>0$, $\tilde{C}_n(\eps, R) >0$ and a positive universal constant $\tilde{\kappa}_n >0$ only depending on the dimension such that
\begin{equation*}
\forall N \geq 1, \ \forall f \in \mathcal{E}_N, \quad \|f\|_{L^2(\rr^n)} \leq \kappa_n(m, R, \gamma, \eps) \Big( \frac{\tilde{\kappa}_n}{\gamma} \Big)^{\tilde{C}_n(\eps, R) N^{1-\frac{\eps}{2}}} \|f\|_{L^2(\omega)},
\end{equation*}
with $\mathcal E_{N}$ being the finite dimensional vector space spanned by the Hermite functions $(\Phi_{\alpha})_{\val \alpha \leq N}$.
\end{theorem}

\medskip

By using the equivalence of norms in finite dimension while taking the parameter $\eps=1$, Theorem~\ref{Spectral} allows to recover the quantitative spectral estimate of Logvinenko-Sereda type (\ref{berk12}) established in~\cite{kkj} (Theorem~2.1), as condition (\ref{thick_rho0}) is then equivalent to the thickness property (\ref{thick1v}). Contrary to the thick case (case $\eps=1$), notice that 
in the case when $0<\eps <1$, the condition (\ref{thick_rho0}) allows control subsets to have holes with diameters tending to infinity.
Theorem~\ref{Spectral} applies for instance with the family of unbounded densities 
$$\rho_\eps(x)= R_\eps \left\langle x \right\rangle^{1-\eps}, \quad x \in \rr^n,$$ 
with $\langle x \rangle=(1+|x|^2)^{\frac{1}{2}}$ and $|\cdot|$ the Euclidean norm on $\rr^n$, when $0 < \eps <1$ and $0< R_\eps \leq \frac{1}{2(1-\eps)}$, as $\rho_\eps$ is then a $\frac{1}{2}$-Lipschitz positive function, see Section~\ref{vsm}. However, the case $\eps =0$ corresponding to a possible linear dependence of the radius is not covered by Theorem~\ref{Spectral}.

The following result shows that the regularity assumptions on the density $\rho$ can be slightly weakened by allowing it to fail being a Lipschitz function, while strengthening on the $\gamma$-thickness condition with respect to $\rho$ by imposing some constraints on the lower bound for the parameter $0<\gamma \leq 1$:

\medskip
 
\begin{corollary}{\label{Spectral2}}
Let $\rho : \rr^n \longrightarrow (0,+\infty)$ be a continuous positive function verifying
 \begin{equation}
\exists 0< \eps \leq 1, \exists 0 < R_{\eps} \leq \frac{1}{2(1-\eps)},  \forall x \in \rr^n, \quad 0<\rho(x) \leq R_{\eps} {\left\langle x\right\rangle}^{1-\eps},
 \end{equation}
with by convention no upper bound condition on $R_\eps>0$ in the case when $\eps=1$.
If $\omega$ is a measurable subset of $\rr^n$ that is $\gamma$-thick with respect to the density $\rho$, that is,
\begin{equation}\label{thick3}
\forall x \in \rr^n, \quad |\omega \cap B(x,\rho(x))| \geq \gamma |B(x,\rho(x))|,
\end{equation}
with $1-\frac{1}{6^n}<\gamma \leq 1$, 
where $B(y,r)$ denotes the Euclidean ball centered at $y \in \rr^n$ with radius $r>0$,
then there exist some positive constants $\kappa_n \big(R_\eps, \gamma, \eps \big)>0$, $\tilde{C}_n (\eps, R_{\eps}) >0$ and a positive universal constant $\tilde{\kappa}_n>0$ only depending on the dimension such that 
\begin{equation*}
\forall N \geq 1, \ \forall f \in \mathcal{E}_N, \quad \|f\|_{L^2(\rr^n)} \leq \kappa_n (R_{\eps}, \gamma, \eps) \Big( \frac{\tilde{\kappa}_n}{\gamma} \Big)^{\tilde{C}_n (\eps, R_{\eps}) N^{1-\frac{\eps}{2}}} \|f\|_{L^2(\omega)}.
\end{equation*}
\end{corollary}

\medskip

The lower bound condition $1-\frac{1}{6^n}<\gamma \leq 1$ can be unexpected. We actually do not know if this assumption is really relevant, or if Corollary~\ref{Spectral2} holds true as well without this technical condition. Let us only mention that this lower bound condition is somehow related with the smallness condition on the positive parameter $0<\eps \leq \eps_0$, with $0<\eps_0 \ll 1$ sufficiently small, in the result of Kovrijkine~\cite{Kovrijkine2} (Theorem~1.1), where is established that a pair $(S, \Sigma)$ is a strong annihilating one when $S$ and $\Sigma$ are measurable subsets satisfying the following $\eps$-thinness condition 
\begin{equation}
\forall x \in \rr^n, \quad |S \cap B(x,\rho_1 (|x|))| \leq \eps |B(x,\rho_1(|x|))|,
\end{equation}
\begin{equation}
\forall x \in \rr^n, \quad |\Sigma \cap B(x,\rho_2 (|x|))| \leq \eps |B(x,\rho_2(|x|))|,
\end{equation}
when $\rho_1, \rho_2 : \rr_+ \longrightarrow (0,+\infty)$ are continous non-increasing functions satisfying 
\begin{equation*}
\exists C_1, C_2>0, \forall t \in \rr_+, \quad \frac{C_2}{\rho_2\Big(\frac{C_1}{\rho_1(t)} \Big)} \geq t,
\end{equation*}
with $0<\eps \leq \eps_0$. Corollary~\ref{Spectral2} is a direct consequence of Theorem~\ref{Spectral} while using the density 
$\rho_\eps(x)=R_{\eps} {\left\langle x\right\rangle}^{1-\eps}$, with $x \in \rr^n$ and $ 0 < R_{\eps} \leq \frac{1}{2(1-\eps)}$, 
together with Lemma~\ref{thick_comparison} in appendix.


\subsection{Null-controllability of hypoelliptic non-selfadjoint quadratic equations}\label{quadratic1} 
This section is devoted to the study of the null-controllability for evolution equations associated to certain classes of non-selfadjoint quadratic operators enjoying some global subelliptic properties. The main result in this section is Corollary~\ref{th2_bis}. This result is a consequence of the new uncertainty principles established in Theorem~\ref{Spectral}, and the abstract observability result given by Theorem~\ref{Meta_thm_AdaptedLRmethod}. It extends to 
any control subset that is thick with respect to an unbounded Lipschitzian density with an almost linear growth at infinity, the result of null-controllability proved by Beauchard, Jaming and the second author in~\cite{kkj} (Theorem~2.2).

\subsubsection{Miscellaneous facts about quadratic differential operators}\label{quadratic}

Quadratic operators are pseudodifferential operators defined in the Weyl quantization
\begin{equation}\label{3}
q^w(x,D_x) f(x) =\frac{1}{(2\pi)^n}\int_{\rr^{2n}}{e^{i(x-y) \cdot \xi}q\Big(\frac{x+y}{2},\xi\Big)f(y)dyd\xi}, 
\end{equation}
by symbols $q(x,\xi)$, with $(x,\xi) \in \rr^{n} \times \rr^n$, $n \geq 1$, which are complex-valued quadratic forms 
\begin{eqnarray*}
q : \rr_x^n \times \rr_{\xi}^n &\rightarrow& \cc\\
 (x,\xi) & \mapsto & q(x,\xi).
\end{eqnarray*}
These operators are actually differential operators with simple and fully explicit expression since the Weyl quantization of the quadratic symbol $x^{\alpha} \xi^{\beta}$, with $(\alpha,\beta) \in \nn^{2n}$, $|\alpha+\beta| = 2$, is given by the differential operator
$$\frac{x^{\alpha}D_x^{\beta}+D_x^{\beta} x^{\alpha}}{2}, \quad D_x=i^{-1}\partial_x.$$
Notice that these operators are non-selfadjoint as soon as their Weyl symbols have a non-zero imaginary part. 
The maximal closed realization of the quadratic operator $q^w(x,D_x)$ on $L^2(\rr^n)$, that is, the operator equipped with the domain
\begin{equation}\label{dom1}
D(q^w)=\big\{f \in L^2(\rr^n) : \ q^w(x,D_x)f \in L^2(\rr^n)\big\},
\end{equation}
where $q^w(x,D_x)f$ is defined in the distribution sense, is known to coincide with the graph closure of its restriction to the Schwartz space~\cite{mehler} (pp.~425-426),
$$q^w(x,D_x) : \mathscr{S}(\rr^n) \rightarrow \mathscr{S}(\rr^n).$$
Classically, to any quadratic form $q : \rr_x^n \times \rr_{\xi}^n \rightarrow \mathbb{C}$ defined on the phase space is associated a matrix $F \in M_{2n}(\CC)$ called its Hamilton map, or its fundamental matrix, which is the unique matrix satisfying the identity
\begin{equation}\label{10}
\forall  (x,\xi) \in \R^{2n},\forall (y,\eta) \in \R^{2n}, \quad q((x,\xi),(y,\eta))=\sigma((x,\xi),F(y,\eta)),
\end{equation}
where $q(\cdot,\cdot)$ is the polarized form associated with the quadratic form $q$, and where $\sigma$ stands for the standard symplectic form
\begin{equation}\label{11}
\sigma((x,\xi),(y,\eta))=\langle \xi, y \rangle -\langle x, \eta\rangle=\sum_{j=1}^n(\xi_j y_j-x_j \eta_j),
\end{equation}
with $x=(x_1,...,x_n)$, $y=(y_1,....,y_n)$, $\xi=(\xi_1,...,\xi_n)$, $\eta=(\eta_1,...,\eta_n) \in \cc^n$. We observe from the definition that 
$$F=\frac{1}{2}\left(\begin{array}{cc}
\nabla_{\xi}\nabla_x q & \nabla_{\xi}^2q  \\
-\nabla_x^2q & -\nabla_{x}\nabla_{\xi} q 
\end{array} \right),$$
where the matrices $\nabla_x^2q=(a_{i,j})_{1 \leq i,j \leq n}$,  $\nabla_{\xi}^2q=(b_{i,j})_{1 \leq i,j \leq n}$, $\nabla_{\xi}\nabla_x q =(c_{i,j})_{1 \leq i,j \leq n}$,
$\nabla_{x}\nabla_{\xi} q=(d_{i,j})_{1 \leq i,j \leq n}$ are defined by the entries
$$a_{i,j}=\partial_{x_i,x_j}^2 q, \quad b_{i,j}=\partial_{\xi_i,\xi_j}^2q, \quad c_{i,j}=\partial_{\xi_i,x_j}^2q, \quad d_{i,j}=\partial_{x_i,\xi_j}^2q.$$
The notion of singular space was introduced in~\cite{kps2} by Hitrik and the second author by pointing out the existence of a particular vector subspace in the phase space $S \subset \rr^{2n}$, which is intrinsically associated with a given quadratic symbol~$q$. This vector subspace  
is defined as the following finite intersection of kernels
\begin{equation}\label{h1bis}
S=\Big( \bigcap_{j=0}^{2n-1}\textrm{Ker}
\big[\textrm{Re }F(\textrm{Im }F)^j \big]\Big)\cap \rr^{2n},
\end{equation}
where $\textrm{Re }F$ and $\textrm{Im }F$ stand respectively for the real and imaginary parts of the Hamilton map $F$ associated with the quadratic symbol $q$,
$$\textrm{Re }F=\frac{1}{2}(F+\overline{F}), \quad \textrm{Im }F=\frac{1}{2i}(F-\overline{F}).$$
As pointed out in \cite{kps2,short,HPSVII,kps11,kps21,rodwahl,viola0}, the notion of singular space plays a basic role in the understanding of the spectral and hypoelliptic properties of the (possibly) non-elliptic quadratic operator $q^w(x,D_x)$, as well as the spectral and pseudospectral properties of certain classes of degenerate doubly characteristic pseudodifferential operators~\cite{kps3,kps4,viola1,viola2}. In particular, the work~\cite{kps2} (Theorem~1.2.2) provides a complete description for the spectrum of any non-elliptic quadratic operator $q^w(x,D_x)$ whose Weyl symbol $q$ has a non-negative real part $\textrm{Re }q \geq 0$, and satisfies a condition of partial ellipticity along its singular space~$S$,
\begin{equation}\label{sm2}
(x,\xi) \in S, \quad q(x,\xi)=0 \Rightarrow (x,\xi)=0. 
\end{equation}
Under these assumptions, the spectrum of the quadratic operator $q^w(x,D_x)$ is shown to be composed of a countable number of eigenvalues with finite algebraic multiplicities and the structure of this spectrum is similar to the one known for elliptic quadratic operators~\cite{sjostrand}. This condition of partial ellipticity is generally weaker than the condition of ellipticity, $S \subsetneq \rr^{2n}$, and allows one to deal with more degenerate situations. An important class of quadratic operators satisfying condition (\ref{sm2}) are those with zero singular spaces $S=\{0\}$. In this case, the condition of partial ellipticity trivially holds.
More specifically, these quadratic operators have been shown in \cite{kps21} (Theorem~1.2.1) to be hypoelliptic and to enjoy global subelliptic estimates of the type
\begin{multline}\label{lol1}
\exists C>0, \forall f \in \mathscr{S}(\rr^n), \\ \|\langle(x,D_x)\rangle^{2(1-\delta)} f\|_{L^2(\rr^n)} \leq C(\|q^w(x,D_x) f\|_{L^2(\rr^n)}+\|f\|_{L^2(\rr^n)}),
\end{multline}
where $\langle(x,D_x)\rangle^{2}=1+|x|^2+|D_x|^2$, with a sharp loss of derivatives $0 \leq \delta<1$ with respect to the elliptic case (case $\delta=0$), which can be explicitly derived from the structure of the singular space.

In this work, we study the class of quadratic operators whose Weyl symbols have non-negative real parts $\textrm{Re }q \geq 0$, and zero singular spaces $S=\{0\}$. 
These quadratic operators are also known~\cite{kps2} (Theorem~1.2.1) to generate strongly continuous contraction semigroups $(e^{-tq^w})_{t \geq 0}$ on $L^2(\rr^n)$, which are smoothing in the Schwartz space for any positive time
$$\forall t>0, \forall f \in L^2(\rr^n), \quad e^{-t q^w}f \in \mathscr{S}(\rr^n).$$
In the recent work~\cite{HPSVII} (Theorem~1.2), these regularizing properties were sharpened and 
these contraction semigroups were shown to be actually smoothing for any positive time in the Gelfand-Shilov space $S_{1/2}^{1/2}(\rr^n)$: $\exists C>0$, $\exists t_0 > 0$, $\forall f \in L^2(\rr^n)$, $\forall \alpha, \beta \in \nn^n$, $\forall 0<t \leq t_0$,
\begin{equation}\label{eq1.10}
\|x^{\alpha}\partial_x^{\beta}(e^{-tq^w}f)\|_{L^{\infty}(\rr^n)} \leq \frac{C^{1+|\alpha|+|\beta|}}{t^{\frac{2k_0+1}{2}(|\alpha|+|\beta|+2n+s)}}(\alpha!)^{1/2}(\beta!)^{1/2}\|f\|_{L^2(\rr^n)},
\end{equation}
where $s$ is a fixed integer verifying $s > n/2$, and where $0 \leq k_0 \leq 2n-1$ is the smallest integer satisfying 
\begin{equation}\label{h1bis2}
\Big( \bigcap_{j=0}^{k_0}\textrm{Ker}\big[\textrm{Re }F(\textrm{Im }F)^j \big]\Big)\cap \rr^{2n}=\{0\}.
\end{equation}
Thanks to this Gelfand-Shilov smoothing effect (\ref{eq1.10}), Beauchard and the second author have established in~\cite{KK1} (Proposition~4.1) that, for any quadratic form $q : \rr_{x,\xi}^{2n} \rightarrow \cc$ with a non-negative real part $\textrm{Re }q \geq 0$ and a zero singular space $S=\{0\}$, the following dissipation estimate holds
\begin{multline}\label{eq6}
\exists C_0>1, \exists t_0>0, \forall t \geq 0, \forall k \geq 0, \forall f \in L^2(\rr^n),\\
\|(1-\pi_k)(e^{-tq^w}f)\|_{L^2(\rr^n)} \leq C_0e^{- \delta(t)k} \|f\|_{L^2(\rr^n)},
\end{multline} 
with 
\begin{equation}\label{eq5}
\delta(t)=\frac{\inf(t,t_0)^{2k_0+1}}{C_0} \geq 0, \qquad t \geq 0, 
\end{equation} 
where $0 \leq k_0 \leq 2n-1$ is the smallest integer satisfying (\ref{h1bis2}), and where $(\pi_k)_{k \geq 0}$ are the Hermite orthogonal projections defined in (\ref{projection}).
Combining these dissipation estimates with the quantitative spectral estimate of Logvinenko-Sereda type (\ref{berk12}) established in~\cite{kkj} (Theorem~2.1), Beauchard, Jaming and the second author have derived from the abstract observability result~\cite{KK1} (Theorem~2.1) the following result of null-controllability~\cite{kkj} (Theorem~2.2):

\medskip

\begin{theorem}[Beauchard, Jaming \& Pravda-Starov]\label{th2}
Let $q : \rr_{x}^{n} \times \rr_{\xi}^n \rightarrow \cc$ be a complex-valued quadratic form with a non-negative real part $\emph{\textrm{Re }}q \geq 0$, and a zero singular space $S=\{0\}$. 
If  $\omega$ is a measurable thick subset of $\mathbb{R}^n$, that is, when condition \emph{(\ref{thick1v})} holds for some $L>0$ and $0<\gamma \leq 1$, then the evolution equation 
$$\left\lbrace \begin{array}{ll}
\partial_tf(t,x) + q^w(x,D_x)f(t,x)=\un_{\omega}(x)u(t,x), \quad &  x \in \mathbb{R}^n, \ t>0,\\
f|_{t=0}=f_0 \in L^2(\rr^n),                                       &  
\end{array}\right.$$
with $q^w(x,D_x)$ being the quadratic differential operator defined by the Weyl quantization of the symbol $q$, is null-controllable from the set $\omega$ in any positive time $T>0$.
\end{theorem}

\medskip

Thanks to the new uncertainty principles established in Theorem~\ref{Spectral}, and the abstract observability result given by Theorem~\ref{Meta_thm_AdaptedLRmethod}, Theorem~\ref{th2} can be generalized to any control subset that is thick with respect to an unbounded Lipschitzian density with an almost linear growth at infinity.

If $\rho : \rr^n \longrightarrow (0,+\infty)$ is a $\frac{1}{2}$-Lipschitz positive function with $\rr^n$ being equipped with the Euclidean norm satisfying that there exist some positive constants $0 < \eps \leq 1$, $m>0$, $R >0$ such that
\begin{equation*}
\forall x \in \rr^n, \quad 0 < m \leq \rho (x) \leq R \left\langle x \right\rangle^{1-\eps}
\end{equation*}
and if $\omega \subset \rr^n$ is a measurable subset that is $\gamma$-thick with respect to the density $\rho$ for some $0<\gamma \leq 1$, that is, when condition (\ref{thick_rho0}) holds, we can apply Theorem~\ref{Meta_thm_AdaptedLRmethod} together with Theorem~\ref{Spectral} for the following choices of parameters:
$\Omega=\mathbb{R}^n$; $A=q^w(x,D_x)$; $0<a=1-\frac{\eps}{2}< b=1$; $t_0>0$ as in (\ref{eq6}); $m_1=2k_0+1$ where $k_0$ is defined in (\ref{h1bis2});
$m_2=0$; any constant $c_1>0$ satisfying
$$\forall k \geq 1,  \quad \kappa_n (m, R, \gamma, \eps) \Big( \frac{\tilde{\kappa}_n}{\gamma} \Big)^{\tilde{C}_n (\eps, R)k^{1-\frac{\eps}{2}}} \leq e^{c_1 k^{1-\frac{\eps}{2}}},$$ 
where the positive constants $\kappa_n (m, R, \gamma, \eps)>0$, $\tilde{C}_n (\eps, R)$, $\tilde{\kappa}_n>0$ are given by Theorem~\ref{Spectral}; $c_1'=c_2'=1$; $c_2= \frac{1}{C_0}>0$ where $C_0>1$ is defined in (\ref{eq6}).
We therefore obtain the following observability estimate in any positive time
\begin{multline*}
\exists C>1, \forall T>0, \forall g \in L^2(\rr^n), \\ \|e^{-Tq^w} g\|_{L^2(\rr^n)}^2 \leq C\exp\Big(\frac{C}{T^{(\frac{2}{\eps}-1)(2k_0+1)}}\Big) \int_0^T \|e^{-tq^w} g\|_{L^2(\omega)}^2 dt.
\end{multline*}
By noticing on one hand that the $L^2(\rr^n)$-adjoint of the quadratic operator $(q^w,D(q^w))$ is the quadratic operator $(\overline{q}^w,D(\overline{q}^w))$, whose Weyl symbol is the complex conjugate of~$q$, and that on the other hand, the symbol $\overline{q}$ is a also a complex-valued quadratic form with a non-negative real part and a zero singular space, the Hilbert Uniqueness Method allows to obtain the following result of null-controllability:

\medskip

\begin{corollary}\label{th2_bis}
Let $q : \rr_{x}^{n} \times \rr_{\xi}^n \rightarrow \cc$ be a complex-valued quadratic form with a non negative real part $\emph{\textrm{Re }}q \geq 0$, and a zero singular space $S=\{0\}$. Let $\rho : \rr^n \longrightarrow (0,+\infty)$ be a $\frac{1}{2}$-Lipschitz positive function with $\rr^n$ being equipped with the Euclidean norm, satisfying that there exist some positive constants $0 < \eps \leq 1$, $m>0$, $R >0$ such that
\begin{equation*}
\forall x \in \rr^n, \quad 0 < m \leq \rho (x) \leq R \left\langle x \right\rangle^{1-\eps}
\end{equation*}
and $\omega$ be a measurable subset of $\mathbb{R}^n$. If $\omega$ is $\gamma$-thick with respect to the density $\rho$, that is, 
$$\exists 0< \gamma \leq 1, \forall x \in \rr^n, \quad |\omega \cap B(x,\rho(x))| \geq \gamma  |B(x,\rho(x))|,$$
where $B(y,r)$ denotes the Euclidean ball centered at $y \in \rr^n$ with radius $r>0$,
then the evolution equation 
$$\left\lbrace \begin{array}{ll}
\partial_tf(t,x) + q^w(x,D_x)f(t,x)=\un_{\omega}(x)u(t,x), \quad &  x \in \mathbb{R}^n, \ t>0, \\
f|_{t=0}=f_0 \in L^2(\rr^n),                                       &  
\end{array}\right.$$
with $q^w(x,D_x)$ being the quadratic differential operator defined by the Weyl quantization of the symbol $q$, is null-controllable from the control subset $\omega$ in any positive time $T>0$.
\end{corollary}

\medskip

\subsection{Null-controllability of evolution equations enjoying Gelfand-Shilov smoothing effects}\label{harmo_oscill}
Given an abstract evolution equation enjoying some Gelfand-Shilov regularizing effects, we aim now at figuring out sufficient geometric conditions on control subsets to ensure its null-controllability in any positive time. 

Let us consider the evolution equation
\begin{equation}\label{PDEfrac}
\left\lbrace \begin{array}{ll}
\partial_tf(t,x) + Af(t,x)=\un_{\omega}(x)u(t,x), \quad &  x \in \mathbb{R}^n, \ t>0, \\
f|_{t=0}=f_0 \in L^2(\rr^n),                                       &  
\end{array}\right.
\end{equation}
associated to $A$ a closed operator on $L^2(\rr^n)$ that is the infinitesimal generator of a strongly continuous contraction semigroup $(e^{-tA})_{t \geq 0}$ on $L^2(\rr^n)$ enjoying some Gelfand-Shilov smoothing effects for any positive time, that is, verifying
\begin{equation}\label{base2}
\forall t>0, \forall u \in L^2(\rr^n), \quad e^{-tA}u \in S^{1/(2s)}_{1/(2s)} (\rr^n),
\end{equation}
with $\frac{1}{2}<s \leq 1$.
We assume more specifically that the contraction semigroup $(e^{-tA})_{t \geq 0}$ enjoys the following quantitative regularizing estimates: 
there exist some constants 
$\frac{1}{2}<s \leq 1$, $C_s > 1$, $0<t_0 \leq 1$, $m_1,m_2 \in \rr$ with $m_1>0$, $m_2 \geq 0$  such that
\begin{multline}\label{GS0}
\forall 0< t \leq t_0, \forall \alpha, \beta  \in \nn^n, \forall g \in L^2(\rr^n),\\
\| x^{\alpha} \partial_x^{\beta}( e^{-t A}g)\| \leq  \frac{C_s^{1+|\alpha|+|\beta|}}{t^{m_1 (|\alpha|+|\beta|)+m_2}} (\alpha !)^{\frac{1}{2s}}(\beta!)^{\frac{1}{2s}}  \|g\|_{L^2(\rr^n)},
\end{multline}
where the above norm $\|\cdot\|$ denotes either the $L^{\infty}(\rr^n)$ norm or the $L^2(\rr^n)$ one. Lemma~\ref{equivv} in appendix shows that if the estimates (\ref{GS0}) hold with the $L^{\infty}(\rr^n)$ norm then they also hold with the $L^2(\rr^n)$ one with the same constants $\frac{1}{2}<s \leq 1$, $0<t_0 \leq 1$, but with different values for the constants $C_s>1$, $m_1>0$, $m_2 \geq 0$. The following result provides sufficient geometric conditions on control subsets related to the index of the symmetric Gelfand-Shilov regularity $\frac{1}{2s}$ to ensure the null-controllability of the adjoint system:

\medskip
 
\begin{theorem}\label{control_PDEfrac}
Let $A$ be a closed operator on $L^2(\rr^n)$ which is the  infinitesimal generator of a strongly continuous contraction semigroup $(e^{-tA})_{t \geq 0}$ on $L^2(\rr^n)$ that satisfies the quantitative smoothing estimates \emph{(\ref{GS0})} for some $\frac{1}{2}< s \leq1$. Let $\rho : \rr^n \longrightarrow (0,+\infty)$ be a $\frac{1}{2}$-Lipschitz positive function
with~$\rr^n$ being equipped with the Euclidean norm,
satisfying that there exist some constants $0 \leq  \delta < 2s-1$, $m>0$, $R>0$ such that
\begin{equation*}
\forall x \in \rr^n, \quad 0<m \leq \rho(x) \leq R{\left\langle x\right\rangle}^{\delta}.
\end{equation*} 
If $\omega$ is a measurable subset of $\rr^n$ which is $\gamma$-thick with respect to the density $\rho$, that is,  
$$\exists 0< \gamma \leq 1, \forall x \in \rr^n, \quad |\omega \cap B(x,\rho(x))| \geq \gamma |B(x,\rho(x))|,$$ 
where $B(y,r)$ denotes the Euclidean ball centered at $y \in \rr^n$ with radius $r>0$,
then the evolution equation associated to the $L^2(\rr^n)$-adjoint operator $A^*$,
 \begin{equation}\label{ucla1}
\left\lbrace \begin{array}{ll}
\partial_tf(t,x) + A^* f(t,x)=\un_{\omega}(x)u(t,x), \quad &  x \in \mathbb{R}^n, \ t>0, \\
f|_{t=0}=f_0 \in L^2(\rr^n),                                       &  
\end{array}\right.
\end{equation}
is null-controllable from the control subset $\omega$ in any positive time $T>0$. 
\end{theorem}

\medskip

As recalled in the previous section, strongly continuous contraction semigroups generated by accretive non-selfadjoint quadratic operators with zero singular spaces enjoy smoothing effects in the Gelfand-Shilov space $S^{1/2}_{1/2}(\rr^n)$. More specifically, Alphonse and Bernier have established in~\cite{P_Alphonse} (Theorem~1.6) that such contraction semigroups $(e^{-tq^{w}})_{t \geq 0}$ on $L^2(\rr^n)$ satisfy the following quantitative regularizing estimates: there exist some constants 
$C > 1$, $0<t_0 \leq 1$ such that
\begin{multline}\label{GS044}
\forall 0< t \leq t_0, \forall k \geq 1, \forall X_1,...,X_k \in \rr^{2n}, \forall g \in L^2(\rr^n),\\
\| L_{X_1}...L_{X_k}( e^{-t q^{w}}g)\|_{L^2(\rr^n)} \leq  \frac{C^{1+k}}{t^{\frac{2k_0+1}{2}k}} \Big(\prod_{j=1}^k |X_j| \Big)(k!)^{\frac{1}{2}} \|g\|_{L^2(\rr^n)}, 
\end{multline}
with $0 \leq k_0 \leq 2n-1$ the smallest integer satisfying (\ref{h1bis2}), where $|X_0|$ is the Euclidean norm of $X_0 \in \rr^{2n}$, and where $L_{X_j}$ is the first order differential operator 
$$L_{X_j}=\langle x_j, x \rangle + \langle \xi_j, \partial_x \rangle, \qquad X_j=(x_j,\xi_j) \in \rr^{2n},$$
with $\langle \cdot,\cdot\rangle$ the Euclidean dot product.
The estimates (\ref{GS044}) imply in particular that for all $0< t \leq t_0$, $\alpha, \beta \in \nn^n$, $g \in L^2(\rr^n)$,
\begin{equation}\label{base3}
 \quad \| x^{\alpha} \partial_x^{\beta}(e^{-t q^{w}}g)\|_{L^2(\rr^n)} \leq \frac{C((2n)^{\frac{1}{2}}C)^{|\alpha|+|\beta|}}{t^{\frac{2k_0+1}{2}(|\alpha|+|\beta|)}}(\alpha !)^{\frac{1}{2}}(\beta!)^{\frac{1}{2}} \|g\|_{L^2(\rr^n)}.
\end{equation}
Indeed, we observe that 
$$x^{\alpha}\partial_x^{\beta}= 
\Big(\prod_{j=1}^n L_{e_j}^{\alpha_j}\Big)\Big(\prod_{k=1}^n L_{\eps_k}^{\beta_k}\Big),\qquad \alpha=(\alpha_1,...,\alpha_n),\ \beta=(\beta_1,...,\beta_n) \in \nn^n,$$ 
where $(e_1,...,e_n,\eps_1,...,\eps_n)$ denotes the canonical basis of $\rr_x^n\times \rr_{\xi}^n$, and that the basic estimate (\ref{base1}) implies that  
$$\forall \alpha, \beta \in \nn^n, \quad (|\alpha|+|\beta|)! \leq 2^{|\alpha|+|\beta|}(|\alpha|)!(|\beta|)! \leq (2n)^{|\alpha|+|\beta|}\alpha!\beta!,$$
since
\begin{equation}\label{fo1}
\frac{(|\alpha|+|\beta|)!}{(|\alpha|)!(|\beta|)!}=\binom{|\alpha|+|\beta|}{|\alpha|} \leq \sum_{k=0}^{|\alpha|+|\beta|}\binom{|\alpha|+|\beta|}{k}=2^{|\alpha|+|\beta|}.
\end{equation}
The strongly continuous contraction semigroup generated by the $L^2(\rr^n)$-adjoint operator $(q^w)^*=(\overline{q})^w$ satisfies the very same quantitative regularizing estimates (\ref{GS044}), since the quadratic symbol $\overline{q}$ has also a non-negative real part with a zero singular space.
Thanks to these smoothing estimates, the result of Corollary~\ref{th2_bis} can therefore be recovered while applying Theorem~\ref{control_PDEfrac}.

As noticed at the end of the proof of Theorem~\ref{control_PDEfrac}, the conclusions of Theorem~\ref{control_PDEfrac} holds true as well when the quantitative regularizing estimates (\ref{GS0}) holding for some $\frac{1}{2}< s \leq1$ are replaced by the following assumption 
\begin{multline}{\label{decay_hermite}}
\exists m_1, m_2>0, \exists C_1, C_2 >0, \exists 0<t_0 \leq 1, \forall 0 < t \leq t_0, \forall g \in L^2(\rr^n),\\ \sum_{\alpha \in \nn^n}e^{\frac{2 t^{m_1}}{C_1} (2 |\alpha|+n)^s} |\langle e^{-t A}g, \Phi_\alpha\rangle_{L^2(\rr^n)}|^2  \leq \frac{C_2^2}{t^{2m_2}} \|g\|_{L^2(\rr^n)}^2,
\end{multline}
with  $(\Phi_{\alpha})_{\alpha \in \nn^n}$ the $L^2(\rr^n)$-Hermite basis.
As an application of this remark, we consider the fractional harmonic operator 
\begin{equation}
\forall u \in D\big(\mathcal{H}^s \big), \quad \mathcal{H}^s u=(-\Delta_x+|x|^2)^su = \sum_{\alpha \in \nn^n} (2 |\alpha|+n)^{s} \langle u,\Phi_\alpha \rangle_{L^2(\rr^n)} \Phi_\alpha,
\end{equation}
with $\frac{1}{2}< s \leq1$, equipped with the domain
\begin{equation}\label{domHs}
D\big(\mathcal{H}^s \big)= \Big\{ u \in L^2(\rr^n): \ \sum_{\alpha \in \nn^n} (2 |\alpha|+n)^{2s} |\langle u,\Phi_\alpha\rangle_{L^2(\rr^n)}|^2 < + \infty \Big\}.
\end{equation}
The fractional harmonic oscillator $\mathcal{H}^s$ is a selfadjoint operator generating a strongly continuous contraction semigroup $(e^{-t \mathcal{H}^s})_{t \geq 0}$ on $L^2(\rr^n)$ explicitly given by 
\begin{equation}\label{SG_frac}
\forall t \geq 0, \forall u \in L^2(\rr^n), \quad e^{-t \mathcal{H}^s} u = \sum_{\alpha \in \nn^n} e^{-t(2 |\alpha|+n)^s} \langle u, \Phi_\alpha \rangle_{L^2(\rr^n)} \Phi_\alpha, 
\end{equation}
see e.g.~\cite{Tucsnak} (Propositions~2.6.2 and~2.6.5). As the assumption (\ref{decay_hermite}) trivially holds for the fractional harmonic oscillator, Theorem~\ref{control_PDEfrac} allows to derive the following result of null-controllability:
 
\medskip 
 
\begin{corollary}\label{corofrac}
Let $\frac{1}{2}< s \leq 1$ and $\rho : \rr^n \longrightarrow (0,+\infty)$ be a $\frac{1}{2}$-Lipschitz positive function 
with $\rr^n$ being equipped with the Euclidean norm,
satisfying that there exist some constants $0\leq  \delta < 2s-1$, $m>0$, $R>0$ such that
\begin{equation*}
\forall x \in \rr^n, \quad 0<m \leq \rho(x) \leq R{\left\langle x\right\rangle}^{\delta}.
\end{equation*} If $\omega$ is a measurable subset of $\rr^n$ which is $\gamma$-thick with respect to the density $\rho$, that is, 
$$\exists 0< \gamma \leq 1, \forall x \in \rr^n, \quad |\omega \cap B(x,\rho(x))| \geq \gamma |B(x,\rho(x))|,$$
where $B(y,r)$ denotes the Euclidean ball centered at $y \in \rr^n$ with radius $r>0$,
then the evolution equation associated to the fractional harmonic oscillator $\mathcal{H}^s=(-\Delta_x+|x|^2)^s$,
 \begin{equation*}
\left\lbrace \begin{array}{ll}
\partial_tf(t,x) + \mathcal{H}^s f(t,x)=\un_{\omega}(x)u(t,x), \quad &  x \in \mathbb{R}^n,\ t>0, \\
f|_{t=0}=f_0 \in L^2(\rr^n),                                       &  
\end{array}\right.
\end{equation*}
is null-controllable from the control subset $\omega$ in any positive time $T>0$. 
\end{corollary}
 
\medskip
 
\subsection{Outline of the work}
Section~\ref{proof_spectral} is devoted to the proof of Theorem~\ref{Spectral}. It is the core of the present work. Theorem~\ref{control_PDEfrac} is then proved in Section~\ref{control_PDEfrac_proof}, whereas the appendix in Section~\ref{appendix} gathers miscellaneous facts about the Gamma function, Hermite functions, slowly varying metrics and the Gelfand-Shilov regularity. Some proofs of technical results as Bernstein type estimates are also given in this appendix.

\section{Proof of Theorem~\ref{Spectral}}\label{proof_spectral}

This section is devoted to the proof of Theorem~\ref{Spectral}.
Let $\rho : \rr^n \longrightarrow (0,+\infty)$ be a $\frac{1}{2}$-Lipschitz positive function with $\rr^n$ equipped with the Euclidean norm,
such that there exist some positive constants $0< \eps \leq 1$, $m>0$, $R>0$ satisfying 
\begin{equation}\label{rho_condi}
\forall x \in \rr^n, \quad 0<m \leq \rho(x) \leq R {\left\langle x \right\rangle}^{1- \eps}.
\end{equation}
Let $\omega$ be a measurable subset of $\rr^n$ which is $\gamma$-thick with respect to the density $\rho$, that is,
\begin{equation}\label{thick_rho}
\exists 0 < \gamma \leq 1, \forall x \in \rr^n, \quad |\omega \cap B(x,\rho(x))| \geq \gamma  |B(x,\rho(x))|=\gamma \rho(x)^n |B(0,1)|,
\end{equation}
where $B(x,r)$ denotes the Euclidean ball centered at $x \in \rr^n$ with radius $r>0$, and where $|A|$ denotes the Lebesgue measure of $A$.
Since $\rho$ is a $\frac{1}{2}$-Lipschitz positive function, Lemma~\ref{slowmet} in appendix shows that the family of norms $(\|\cdot\|_x)_{x \in \rr^n}$ given by
\begin{equation}
\forall x \in \rr^n, \forall y \in \rr^n, \quad \|y\|_x=\frac{\|y\|}{\rho(x)},
\end{equation}
where $\|\cdot\|$ denotes the Euclidean norm in $\rr^n$, defines a slowly varying metric on $\rr^n$.

\subsection{Step 1. Bad and good balls}  By using Theorem~\ref{slowmetric} in appendix, we can find a sequence $(x_k)_{k \geq 0}$ in $\rr^n$ such that 
\begin{equation}\label{recov}
\exists N_0 \in \nn, \forall (i_1, ..., i_{N_0+1}) \in \nn^{N_0+1} \textrm{ with } i_k \neq i_l \textrm{ if }1 \leq k \neq l \leq N_0+1, \quad \bigcap \limits_{k=1}^{N_0+1} {B_{i_k}}=\emptyset
\end{equation}
and
\begin{equation}\label{recov1}
 \rr^n=\bigcup_{k=0}^{+\infty} {B_k},
\end{equation}
where 
\begin{equation}\label{asdf1}
B_k=\{y \in \rr^n:\ \|y-x_k\|_{x_k} <1\}=\{y \in \rr^n:\ \|y-x_k\| <\rho(x_k)\}=B(x_k,\rho(x_k)). 
\end{equation}
Let us notice from Theorem~\ref{slowmetric} that the non-negative integer $N_0$ only depends on the dimension $n$ and the constant $C \geq 1$ appearing in slowness condition (\ref{equiv}) which can be taken equal to $C=2$ here, as $\rho$ is a $\frac{1}{2}$-Lipschitz function. The integer $N_0=N_0(n)$ is therefore independent on the function $\rho$ and depends only on the dimension $n \geq 1$.  
It follows from (\ref{recov}) and (\ref{recov1}) that 
\begin{equation}\label{asdf2}
\forall x \in \rr^n, \quad 1 \leq \sum \limits_{k=0}^{+\infty} \mathbbm{1}_{B_k} (x) \leq N_0,
\end{equation}
where $\mathbbm{1}_{B_k}$ denotes the characteristic function of $B_k$.
We deduce from (\ref{asdf2}) that for all $g \in L^2(\rr^n)$,
\begin{equation}
\|g\|_{L^2(\rr^n)}^2 = \int_{\rr^n}|g(x)|^2dx \leq \sum_{k=0}^{+\infty}\int_{B_k}|g(x)|^2dx \leq N_0 \|g\|_{L^2(\rr^n)}^2.
\end{equation}
Let $N \in \nn$ be a non-negative integer and $f \in \mathcal E_{N} \setminus \{0\}$, with $\mathcal E_{N}$ being the finite dimensional vector space spanned by the Hermite functions $(\Phi_{\alpha})_{\val \alpha \leq N}$ defined in (\ref{berk1}).
Let $0<\delta \leq 1$ be a positive constant to be chosen later on.
We divide the family of balls $(B_k)_{k \geq 0}$ into families of good and bad balls. A ball $B_k$, with $k \in \nn$, is said to be good if it satisfies 
\begin{multline}\label{good}
\forall \big(\beta,\tilde{\beta} \big) \in \nn^n \times \nn^n,\ |\tilde{\beta}| \leq n, \\  \int_{B_k}|\left\langle x\right\rangle^{(1-\eps)(\left|\beta\right|+n)}\partial_x^{\beta+\tilde{\beta}}f(x)|^2dx \leq 4^n \big(2(2^n N_0+1)\big)^{\left| \beta \right|+1} M_{\beta,\tilde{\beta},N}(\delta)^2 \int_{B_k}|f(x)|^2dx,
\end{multline}
where the positive constants $M_{\beta,\tilde{\beta},N}(\delta)>0$ also depend on the fixed positive parameter $0<\eps\leq 1$ and the dimension $n \geq 1$, and are defined by
\begin{multline}\label{asdf4}
M_{\beta,\tilde{\beta},N}(\delta) \\ =\tilde{K}_{\eps,\delta} K_\eps^{(2-\eps)|\beta|+(1-\eps)n+|\tilde{\beta}|}\delta^{|\beta|+|\tilde{\beta}|} \big( n+1 \big)^{\frac{(1-\eps)(n+|\beta|)}{2}} \Gamma\Big(|\beta|+\frac{(1-\eps)n+|\tilde{\beta}|}{2-\eps}+3\Big) e^{\frac{N^{1-\frac{\eps}{2}}}{\delta^{2-\eps}}},
\end{multline}
with the constants $\tilde{K}_{\eps,\delta}>1$ and $K_\eps>1$ defined in Proposition~\ref{prop1}.
On the other hand, a ball $B_k$, with $k \in \nn$, which is not good, is said to be bad, that is, when
\begin{multline}\label{bad}
\exists \big(\beta,\tilde{\beta} \big) \in \nn^n \times \nn^n, \ |\tilde{\beta}| \leq n, \\  \int_{B_k}|\left\langle x\right\rangle^{(1-\eps)(\left|\beta\right|+n)}\partial_x^{\beta+\tilde{\beta}}f(x)|^2dx >4^n \big(2(2^n N_0+1)\big)^{\left| \beta \right|+1} M_{\beta,\tilde{\beta},N}(\delta)^2 \int_{B_k}|f(x)|^2dx.
\end{multline}
If $B_k$ is a bad ball, it follows from (\ref{bad}) that there exists $\big(\beta_0,\tilde{\beta}_0 \big) \in \nn^n \times \nn^n$, $|\tilde{\beta}_0|\leq n$ such that 
\begin{align}\label{gh05}
& \ \int_{B_k}|f(x)|^2dx \\ \notag
\leq & \  
\frac{1}{4^n\big(2(2^n N_0+1)\big)^{\left|\beta_0\right|+1}M_{\beta_0,\tilde{\beta}_0,N}(\delta)^2}\int_{B_k}|\left\langle x\right\rangle^{(1-\eps)(|\beta_0|+n)}\partial_x^{\beta_0+\tilde{\beta}_0}f(x)|^2dx \\ \notag
\leq & \  \sum_{\substack{\beta \in \nn^n, \\ \tilde{\beta} \in \nn^n,\ |\tilde{\beta}| \leq n}} \frac{1}{4^n\big(2(2^n N_0+1)\big)^{\left|\beta\right|+1}M_{\beta, \tilde{\beta},N}(\delta)^2}\int_{B_k}|\left\langle x\right\rangle^{(1-\eps)(|\beta|+n)}\partial_x^{\beta+\tilde{\beta}}f(x)|^2dx.
\end{align}
By summing over all the bad balls and by using from (\ref{recov}) that
\begin{equation}
\forall x \in \rr^n, \quad \mathbbm{1}_{\bigcup_{\textrm{bad balls}} B_k} \leq \sum_{\textrm{bad balls}} \mathbbm{1}_{B_k} \leq N_0 \mathbbm{1}_{\bigcup_{\textrm{bad balls}} B_k},
\end{equation} we deduce from (\ref{gh05}) and the Fubini-Tonelli theorem that 
\begin{align}\label{gh6}
& \ \int_{\bigcup_{\textrm{bad balls}} B_k}|f(x)|^2dx \leq \sum_{\textrm{bad balls}}\int_{B_k}|f(x)|^2dx\\ \notag
\leq & \ \sum_{\substack{\beta \in \nn^n \\ \tilde{\beta} \in \nn^n,\ |\tilde{\beta}| \leq n}} \frac{N_0}{4^n \big(2(2^n N_0+1)\big)^{\left|\beta\right|+1}M_{\beta,\tilde{\beta},N}(\delta)^2} \int_{\bigcup_{\textrm{bad balls}}  B_k} \hspace{-8mm}  |{\left\langle x\right\rangle}^{(1-\eps)(\left|\beta\right|+n)} \partial_x^{\beta+\tilde{\beta}} f(x)|^2dx\\ \notag
\leq & \ \sum_{\substack{\beta \in \nn^n \\ \tilde{\beta} \in \nn^n,\ |\tilde{\beta}| \leq n}} \frac{N_0}{4^n \left(2(2^n N_0+1)\right)^{\left|\beta\right|+1}M_{\beta,\tilde{\beta},N}(\delta)^2} \int_{\rr^n}|\left\langle x\right\rangle^{(1-\eps)(\left|\beta\right|+n)}\partial_x^{\beta+\tilde{\beta}} f(x)|^2dx.
\end{align}
By using that the number of solutions to the equation $\beta_1+...+\beta_{n}=k$, with $k \geq 0$, $n \geq 1$ and unknown $\beta=(\beta_1,...,\beta_n) \in \nn^{n}$, is given by $\binom{k+n-1}{k}$, we obtain from the Bernstein type estimates in Proposition~\ref{prop1}, (\ref{asdf4}) and (\ref{gh6}) that 
\begin{align}\label{gh6y}
& \ \int_{\bigcup_{\textrm{bad balls}} B_k}|f(x)|^2dx \leq \Big(\sum_{\substack{\beta \in \nn^n, \\ \tilde{\beta} \in \nn^n,\ |\tilde{\beta}| \leq n}} \frac{N_0}{4^n \big(2(2^{n}N_0+1)\big)^{|\beta|+1}}\Big)
\|f\|_{L^2(\rr^n)}^2\\ \notag
= & \ \Big(\sum_{\beta \in \nn^n}  \frac{N_0}{\big(2(2^{n}N_0+1)\big)^{|\beta|+1}}\Big) \Big(\sum_{\tilde{\beta} \in \nn^n,\ |\tilde{\beta}| \leq n} \frac{1}{4^n} \Big)
\|f\|_{L^2(\rr^n)}^2\\ \notag
=& \ \Big(\sum_{k=0}^{+\infty}\binom{k+n-1}{k} \frac{N_0}{2^{k+1}(2^{n}N_0+1)^{k+1}}\Big) \Big(\sum_{j=0}^{n}\frac{1}{4^n}\binom{j+n-1}{j}\Big)
\|f\|_{L^2(\rr^n)}^2\\ \notag
\leq & \  2^{n-2}\Big(\sum_{k=1}^{+\infty}\frac{N_0}{(2^{n}N_0+1)^{k}}\Big) \Big(\sum_{j=0}^{n}\frac{2^{j+n-1}}{4^n}\Big)
\|f\|_{L^2(\rr^n)}^2 \leq \frac{1}{4}\|f\|_{L^2(\rr^n)}^2,
\end{align}
since 
\begin{equation}\label{gh45}
\binom{k+n-1}{k} \leq \sum_{j=0}^{k+n-1}\binom{k+n-1}{j}=2^{k+n-1}.
\end{equation}
Recalling from (\ref{recov1}) that 
$$ 1 \leq \mathbbm{1}_{\bigcup_{\textrm{bad balls}}B_k}+ \mathbbm{1}_{\bigcup_{\textrm{good balls}}B_k},$$ 
we notice that  
\begin{equation}\label{asdf5}
\|f\|_{L^2(\rr^n)}^2 \leq \int_{\bigcup_{\textrm{good balls}} B_k}|f(x)|^2dx+ \int_{\bigcup_{\textrm{bad balls}} B_k}|f(x)|^2dx.
\end{equation}
It follows from (\ref{gh6y}) and (\ref{asdf5}) that 
\begin{equation}\label{gh7}
\|f\|_{L^2(\rr^n)}^2 \leq \frac{4}{3} \int_{\bigcup_{\textrm{good balls}} B_k}|f(x)|^2dx.
\end{equation}

\subsection{Step 2. Properties on good balls}
As the ball $B(0,1)$ is an Euclidean ball, the Sobolev embedding 
$$W^{n,2}(B(0,1)) \xhookrightarrow{} L^{\infty}(B(0,1)),$$
see e.g.~\cite{adams} (Theorem~4.12), implies that there exists a positive constant $C_{n}>0$ depending only the dimension $n \geq 1$ such that 
\begin{equation}\label{sobolev}
\forall u \in W^{n,2}(B(0,1)), \quad 
\|u\|_{L^{\infty}(B(0,1))} \leq C_{n} \|u\|_{W^{n,2}(B(0,1))}.
\end{equation}
By translation invariance and homogeneity of the Lebesgue measure, it follows from (\ref{rho_condi}), (\ref{asdf1}) and (\ref{sobolev}) that for all $u \in {W^{n,2}(B_k)}$,
\begin{multline*}
\|u\|^2_{L^{\infty}(B_k)}=\|x \mapsto u(x_k+x \rho(x_k))\|^2_{L^{\infty}(B(0,1))}
 \leq C_{n}^2 \|x \mapsto u(x_k+x \rho(x_k))\|^2_{W^{n,2}(B(0,1))} \\
 =C_{n}^2 \sum_{\substack{\alpha \in \nn^n, \\ |\alpha| \leq n}} \int_{B_k} \rho(x_k)^{2|\alpha|-n} |\partial^{\alpha}_x u(x)|^2 dx
 =C_{n}^2 \sum_{\substack{\alpha \in \nn^n, \\ |\alpha| \leq n}} \int_{B_k}m^{2|\alpha|-n} \Big( \frac{\rho(x_k)}{m} \Big)^{2|\alpha|-n} |\partial^{\alpha}_x u(x)|^2 dx\end{multline*}
and
\begin{multline}{\label{se1}}
\|u\|^2_{L^{\infty}(B_k)}
 \leq C_{n}^2 \max(m,m^{-1})^n \sum_{\substack{\alpha \in \nn^n, \\ |\alpha| \leq n}} \int_{B_k} \Big( \frac{\rho(x_k)}{m} \Big)^{n} |\partial^{\alpha}_x u(x)|^2 dx\\
 = C_{n}^2 \max(1,m^{-1})^{2n} \rho(x_k)^{n} \sum_{\substack{\alpha \in \nn^n, \\ |\alpha| \leq n}} \int_{B_k} |\partial^{\alpha}_x u(x)|^2 dx.
\end{multline}
We deduce from (\ref{se1}) that for all $u \in {W^{n,2}(B_k)}$,
\begin{equation}{\label{se2}}
\|u\|_{L^{\infty}(B_k)} \leq C_{n} \max(1,m^{-1})^{n} \rho(x_k)^{\frac{n}{2}} \|u\|_{W^{n,2}(B_k)}.
\end{equation}
Let $B_k$ be a good ball. By using the fact that the mapping $\rho$ is a $\frac{1}{2}$-Lipschitz positive function, we notice that
\begin{equation}{\label{equi}}
\forall x \in B_k=B(x_k,\rho(x_k)), \quad 0 < \rho(x_k) \leq 2 \rho(x).
\end{equation}
We deduce from (\ref{se2}) and (\ref{equi}) that for all $\beta \in \nn^n$ and $k \in \nn$ such that $B_k$ is a good ball
\begin{align}\label{gh30}
& \  \rho(x_k)^{|\beta|+ \frac{n}{2}}\|\partial_x^{\beta}f\|_{L^{\infty}(B_k)} \\ \notag
\leq & \ C_n \max(1,m^{-1})^{n} \rho(x_k)^{|\beta|+ n}\Big(\sum_{\substack{\tilde{\beta} \in \nn^n, \ |\tilde{\beta}| \leq n}}\|\partial_x^{\beta+\tilde{\beta}}f\|^2_{L^{2}(B_k)}\Big)^{\frac{1}{2}}\\ \notag
=  & \ C_n \max(1,m^{-1})^{n}  \Big(\sum_{\substack{\tilde{\beta} \in \nn^n, \ |\tilde{\beta}| \leq n}}\| \rho(x_k)^{|\beta|+ n}\partial_x^{\beta+\tilde{\beta}}f\|^2_{L^{2}(B_k)}\Big)^{\frac{1}{2}} \\ \notag
\leq & \ C_n \max(1,m^{-1})^{n} 2^{|\beta|+n} \Big(\sum_{\substack{\tilde{\beta} \in \nn^n, \ |\tilde{\beta}| \leq n}}\| \rho(x)^{|\beta|+ n}\partial_x^{\beta+\tilde{\beta}}f\|^2_{L^{2}(B_k)}\Big)^{\frac{1}{2}}.
\end{align}
By using (\ref{rho_condi}) and the definition of good balls (\ref{good}), it follows from (\ref{gh30}) that for all $\beta \in \nn^n$ and $k \in \nn$ such that $B_k$ is a good ball
\begin{align}\label{asdf7}
 & \ \rho(x_k)^{|\beta|+ \frac{n}{2}}\|\partial_x^{\beta}f\|_{L^{\infty}(B_k)} \\ \notag
\leq & \ C_n \max(1,m^{-1})^{n} \big(2R\big)^{|\beta|+n} \Big(\sum_{\substack{\tilde{\beta} \in \nn^n, \\ |\tilde{\beta}| \leq n}}\| {\left\langle x \right\rangle}^{(1-\eps)(|\beta|+ n)}\partial_x^{\beta+\tilde{\beta}}f\|^2_{L^{2}(B_k)}\Big)^{\frac{1}{2}} \\
\notag
\leq & \ C_n \max(1,m^{-1})^{n} \big(2R\big)^{|\beta|+n} 2^n \sqrt{2(2^nN_0+1)}^{|\beta|+1} \Big(\sum_{\substack{\tilde{\beta} \in \nn^n, \\ |\tilde{\beta}| \leq n}}
M_{\beta,\tilde{\beta},N}(\delta)^2 \Big)^{\frac{1}{2}}\|f\|_{L^2(B_k)}.
\end{align}

\noindent By using the fact that the Gamma function is increasing on $[2,+\infty)$ (see Section~\ref{miscgamma}), we obtain from (\ref{asdf4}) that 
for all $\beta \in \nn^n$, $\tilde{\beta} \in \nn^n$, $|\tilde{\beta}| \leq n$, $0<\delta \leq 1$,
\begin{equation}{\label{gh35}}
\ M_{\beta,\tilde{\beta},N}(\delta) \leq \tilde{K}_{\eps,\delta} K_{\eps}^{(2-\eps)|\beta|+n(1-\eps)+|\tilde{\beta}|}\delta^{|\beta|+|\tilde{\beta}|} (n+1)^{\frac{(1-\eps)(|\beta|+n)}{2}} \Gamma \big(|\beta|+n+3 \big) e^{\frac{N^{1-\frac{\eps}{2}}}{\delta^{2-\eps}}}.
\end{equation}
Recalling that $K_{\eps}>1$ and $0<\delta \leq 1$, it follows from (\ref{asdf7}) and (\ref{gh35}) that for all $\beta \in \nn^n$ and $k \in \nn$ such that $B_k$ is a good ball
\begin{multline}{\label{gh36}}
\rho(x_k)^{|\beta|+ \frac{n}{2}}\|\partial_x^{\beta}f\|_{L^{\infty}(B_k)} \\
\leq C_n(\delta, \eps, m, R) \big(\delta \tilde{C}_n(\eps, R) \big)^{|\beta|} \Gamma \big( |\beta|+n+3 \big) e^{\frac{N^{1-\frac{\eps}{2}}}{\delta^{2-\eps}}} \|f\|_{L^2(B_k)},
\end{multline}
with 
\begin{multline}\label{gh50}
C_n(\delta, \eps, m, R)\\ =\tilde{K}_{\eps,\delta} C_n \max(1,m^{-1})^{n} (4R)^n  K_{\eps}^{(2-\eps)n} (n+1)^{(1-\eps)\frac{n}{2}} \sqrt{2(2^n N_0+1)}>0
\end{multline}
and 
\begin{equation}
\tilde{C}_n(\eps, R)=2 R\sqrt{2(2^nN_0+1)}  K_{\eps}^{2-\eps} (n+1)^{\frac{1-\eps}{2}}>0.
\end{equation}
Let $B_k$ be a good ball. Recalling that $f$ is a finite combination of Hermite functions, we deduce from the continuity of the function $f$ and the compactness of $\overline{B_k}$ that there exists $y_{k} \in \overline{B_k}$ such that 
\begin{equation}\label{gh8}
\|f\|_{L^{\infty}(B_k)}=|f(y_{k})|.
\end{equation}
By using spherical coordinates centered at $y_{k} \in \overline{B_k}$ and the fact that the Euclidean diameter of the ball $B_k=B(x_k,\rho(x_k))$ is $2\rho(x_k)$, we observe that 

\begin{align}\label{gh9}
|\omega \cap B_k|=& \ \int_0^{+\infty}\Big(\int_{\mathbb{S}^{n-1}}\un_{\omega \cap B_k}(y_k+r \sigma)d\sigma\Big)r^{n-1}dr\\ \notag
= & \ \int_0^{2\rho(x_k)}\Big(\int_{\mathbb{S}^{n-1}}\un_{\omega \cap B_k}(y_k+r \sigma)d\sigma\Big)r^{n-1}dr\\ \notag
= & \ \big(2\rho(x_k)\big)^n \int_0^{1}\Big(\int_{\mathbb{S}^{n-1}}\un_{\omega \cap B_k}(y_k+2\rho(x_k)r \sigma)d\sigma\Big)r^{n-1}dr,
\end{align}

\noindent where $\un_{\omega \cap B_k}$ denotes the characteristic function of the measurable set $\omega \cap B_k$.
By using the Fubini-Tonelli theorem, we deduce from (\ref{gh9}) that 

\begin{multline}\label{gh10}
|\omega \cap B_k| \leq \big(2\rho(x_k)\big)^n \int_0^{1}\Big(\int_{\mathbb{S}^{n-1}}\un_{\omega \cap B_k}(y_k+2\rho(x_k)r \sigma)d\sigma\Big)dr\\
=\big(2\rho(x_k)\big)^n\int_{\mathbb{S}^{n-1}}\Big(\int_0^{1}\un_{\omega \cap B_k}(y_k+2\rho(x_k)r \sigma)dr\Big)d\sigma\\
=\big(2\rho(x_k)\big)^n\int_{\mathbb{S}^{n-1}}\Big(\int_0^{1}\un_{I_{\sigma}}(r)dr\Big)d\sigma=\big(2\rho(x_k)\big)^n\int_{\mathbb{S}^{n-1}}|I_{\sigma}|d\sigma,
\end{multline}

\noindent where

\begin{equation}\label{gh11}
I_{\sigma}=\{r \in [0,1] :\ y_k+2\rho(x_k)r \sigma \in \omega \cap B_k\}.
\end{equation}
The estimate (\ref{gh10}) implies that there exists $\sigma_0(k) \in \mathbb{S}^{n-1}$ such that 
\begin{equation}\label{gh12}
|\omega \cap B_k| \leq \big(2\rho(x_k)\big)^n|\mathbb{S}^{n-1}| |I_{\sigma_0(k)}|.
\end{equation}
Recalling that $B_k=B\big(x_k, \rho(x_k) \big)$ and using the property (\ref{thick_rho}), it follows from (\ref{gh12}) that 
\begin{equation}\label{gh13}
0< \frac{\gamma |B(0,1)|}{2^n|\mathbb{S}^{n-1}|} \leq \frac{|\omega \cap B(x_k,\rho(x_k))|}{(2\rho(x_k))^n|\mathbb{S}^{n-1}|} \leq |I_{\sigma_0(k)}| \leq 1.
\end{equation}

\subsection{Step 3. Recovery of the $L^2(\rr^n)$-norm}
Let $B_k$ be a good ball. We first notice that $\|f\|_{L^2(B_k)} \neq 0$, since $f$ is a non-zero entire function. We consider the entire function
\begin{equation}\label{gh13b}
\forall z \in \cc, \quad \phi(z)=\left|B_k\right|^{\frac{1}{2}}\frac{f(y_k+2\rho(x_k) z \sigma_0(k))}{\|f\|_{L^2(B_k)}},
\end{equation}
where $y_k$ and $\sigma_0(k)$ are defined in (\ref{gh8}) and (\ref{gh12}).
We observe from (\ref{gh8}) that 
$$|\phi(0)|=\left|B_k\right|^{\frac{1}{2}}\frac{|f(y_k)|}{\|f\|_{L^2(B_k)}}=\left|B_k\right|^{\frac{1}{2}} \frac{\|f\|_{L^{\infty}(B_k)}}{{\|f\|_{L^2(B_k)}}} \geq 1.$$
Instrumental in the proof of Theorem~\ref{Spectral} is the following lemma proved by Kovrijkine in~\cite{Kovrijkine} (Lemma~1):

\medskip

\begin{lemma}\label{lem_kov}
Let $I \subset \rr$ be an interval of length $1$ such that $0 \in I$ and $E \subset I$ be a subset of positive measure $|E|>0$. There exists a positive constant $C>1$ such that for all analytic function $\Phi$ on the open ball $B_{\cc}(0,5)$ centered in zero with radius~$5$ such that $|\Phi(0)| \geq 1$, then 
$$\sup_{x \in I}|\Phi(x)| \leq \Big(\frac{C}{|E|}\Big)^{\frac{\ln M}{\ln 2}}\sup_{x \in E}|\Phi(x)| ,$$
with $M=\sup_{|z| \leq 4}|\Phi(z)| \geq 1$.  
\end{lemma}

\medskip

Applying Lemma~\ref{lem_kov} with $I=[0,1]$, $E=I_{\sigma_0(k)} \subset [0,1]$ verifying $|E|=|I_{\sigma_0(k)}|>0$ according to (\ref{gh13}), and the analytic function $\Phi=\phi$ defined in (\ref{gh13b}) satisfying $|\phi(0)| \geq 1$, we obtain that
\begin{multline}\label{gh20}
\left|B_k\right|^{\frac{1}{2}}\frac{\sup_{x \in [0,1]}|f(y_k+2\rho(x_k) x \sigma_0(k))|}{\|f\|_{L^2(B_k)}} \\ \leq \Big(\frac{C}{|I_{\sigma_0(k)}|}\Big)^{\frac{\ln M}{\ln 2}}\left|B_k\right|^{\frac{1}{2}}\frac{\sup_{x \in I_{\sigma_0(k)}}|f(y_k+2\rho(x_k) x \sigma_0(k))|}{\|f\|_{L^2(B_k)}},
\end{multline}
with
\begin{equation}\label{gh22}
1 \leq M=\left|B_k\right|^{\frac{1}{2}}\frac{\sup_{|z| \leq 4}|f(y_k+2\rho(x_k) z \sigma_0(k))|}{\|f\|_{L^2(B_k)}}.
\end{equation}
It follows from (\ref{gh13}) and (\ref{gh20}) that 
\begin{multline}\label{gh21}
\sup_{x \in [0,1]}|f(y_k+2\rho(x_k) x \sigma_0(k))| \leq \Big(\frac{2^n C |\mathbb{S}^{n-1}|}{\gamma |B(0,1)|}\Big)^{\frac{\ln M}{\ln 2}}\sup_{x \in I_{\sigma_0(k)}}|f(y_k+2\rho(x_k) x \sigma_0(k))|\\
\leq M^{\frac{1}{\ln 2}\ln(\frac{2^n C |\mathbb{S}^{n-1}|}{\gamma |B(0,1)|})}
\sup_{x \in I_{\sigma_0(k)}}|f(y_k+2\rho(x_k) x \sigma_0(k))|.
\end{multline}
According to (\ref{gh11}), we notice that 
\begin{equation}\label{gh23}
\sup_{x \in I_{\sigma_0(k)}}|f(y_k+2\rho(x_k) x \sigma_0(k))| \leq \|f\|_{L^{\infty}(\omega \cap B_k)}.
\end{equation}
On the other hand, we deduce from (\ref{gh8}) that 
\begin{equation}\label{gh24}
\|f\|_{L^{\infty}(B_k)}= |f(y_k)| \leq \sup_{x \in [0,1]}|f(y_k+2\rho(x_k) x \sigma_0(k))|.
\end{equation}
It follows from (\ref{gh21}), (\ref{gh23}) and (\ref{gh24}) that 
\begin{equation}\label{gh25}
\|f\|_{L^{\infty}(B_k)}\leq M^{\frac{1}{\ln 2}\ln(\frac{2^n C |\mathbb{S}^{n-1}|}{\gamma |B(0,1)|})}\|f\|_{L^{\infty}(\omega \cap B_k)}.
\end{equation}
By using the analyticity of the entire function $f$, we observe that 
\begin{equation}\label{gh31}
\forall z \in \cc, \quad f(y_k+2\rho(x_k)z \sigma_0(k))
=\sum_{\beta \in \nn^n}\frac{(\partial_x^{\beta}f)(y_k)}{\beta!}\sigma_0(k)^{\beta} \big(2\rho(x_k) \big)^{|\beta|}z^{|\beta|}.
\end{equation}
By using that $B_k=B\big(x_k, \rho(x_k) \big)$ is a good ball, $y_k \in \overline{B_k}$ and the continuity of the functions $\partial_x^{\beta}f$, we deduce from (\ref{gh36}) and (\ref{gh31}) that for all $|z| \leq 4$,
\begin{multline}\label{gh32}
|B_k|^{\frac{1}{2}}|f(y_k+2\rho(x_k)z \sigma_0(k))|= \rho(x_k)^{\frac{n}{2}} |B(0,1)|^{\frac{1}{2}} |f(y_k+2\rho(x_k)z \sigma_0(k))| \\
\leq |B(0,1)|^{\frac{1}{2}} \sum_{\beta \in \nn^n} \rho(x_k)^{|\beta|+\frac{n}{2}}\frac{|(\partial_x^{\beta}f)(y_k)|}{\beta!} 8^{|\beta|}
\leq |B(0,1)|^{\frac{1}{2}}\sum_{\beta \in \nn^n} \rho(x_k)^{|\beta|+\frac{n}{2}}\frac{\|\partial_x^{\beta}f\|_{L^{\infty}(B_k)}}{\beta!} 8^{|\beta|} \\
\leq |B(0,1)|^{\frac{1}{2}}C_n(\delta,\eps, m, R) e^{\frac{N^{1-\frac{\eps}{2}}}{\delta^{2-\eps}}}\Big(\sum_{\beta \in \nn^n}\frac{\Gamma \big(|\beta|+n+3 \big)}{\beta!}\big(8\delta \tilde{C}_n(\eps, R) \big)^{|\beta|} \Big)\|f\|_{L^2(B_k)}.
\end{multline}
We recall the following estimate
\begin{equation}\label{base1}
\forall \beta \in \nn^n, \quad |\beta|! \leq n^{|\beta|} \beta !,
\end{equation}
which is obtained by using the Newton formula, see formula (0.3.3) in~\cite{rodino}.
By using anew that the number of solutions to the equation $\beta_1+...+\beta_{n}=k$, with $k \geq 0$, $n \geq 1$ and unknown $\beta=(\beta_1,...,\beta_n) \in \nn^{n}$, is given by $\binom{k+n-1}{k}$, and that
$$\quad \Gamma \big( |\beta| +n+3 \big)= (|\beta|+n+2)!, \quad (|\beta|+n+2)^{n+2} \leq (n+2)! e^{|\beta|+n+2},$$ 
according to (\ref{factoriel}) and (\ref{rod1}),
we notice from (\ref{gh45}) that 
\begin{multline}\label{gh49}
\sum_{\beta \in \nn^n}\frac{\Gamma \big(|\beta|+n+3 \big)}{\beta!}\big(8 \delta \tilde{C}_n(\eps, R) \big)^{|\beta|} = \sum_{\beta \in \nn^n}\frac{(|\beta|+n+2)!}{\beta!}\big(8 \delta \tilde{C}_n(\eps, R)  \big)^{|\beta|} \\
\leq \sum_{\beta \in \nn^n} (|\beta|+n+2)^{n+2}\frac{|\beta|!}{\beta!}\big(8 \delta \tilde{C}_n(\eps, R) \big)^{|\beta|}
\leq e^{n+2} (n+2)! \sum_{\beta \in \nn^n}\big(8 \delta n e \tilde{C}_n(\eps, R) \big)^{|\beta|}\\
=e^{n+2} (n+2)! \sum_{k=0}^{+\infty}\binom{k+n-1}{k}\big(8\delta n e \tilde{C}_n(\eps, R) \big)^{k}
\leq e^{n+2}(n+2)! 2^{n-1}\sum_{k=0}^{+\infty}\big(16 \delta n e \tilde{C}_n(\eps, R) \big)^{k}.
\end{multline}
We can now make the following choice for the positive parameter $0<\delta \leq 1$, which is fixed from now and taken to be equal to
\begin{equation}\label{gh51}
0<\delta=\delta_{n,\eps,R}= \min\Big(1,\frac{1}{32 n e \tilde{C}_n(\eps, R)}\Big) \leq 1.
\end{equation}
Setting $D_n(\eps, m, R)=C_n\big(\delta_{n,\eps,R},\eps, m, R\big)>0$, it follows from (\ref{gh22}), (\ref{gh32}), (\ref{gh49}) and (\ref{gh51}) that
\begin{equation}\label{gh52}
1 \leq M \leq |B(0,1)|^{\frac{1}{2}}(n+2)! D_n(\eps, m, R) e^{n+2} 2^n e^{\delta_{n,\eps,R}^{\eps-2} N^{1-\frac{\eps}{2}}}.
\end{equation}
We notice from (\ref{gh13}) that  
\begin{equation}\label{gh52s}
\frac{2^n C|\mathbb{S}^{n-1}|}{\gamma |B(0,1)|} > 1,
\end{equation}
since the positive constant given by Lemma~\ref{lem_kov} satisfies $C>1$. With this choice, we deduce from (\ref{gh25}) and (\ref{gh52}) that 
\begin{equation}\label{gh53}
\|f\|_{L^{\infty}(B_k)}\leq 
\Big(\frac{2^n C|\mathbb{S}^{n-1}|}{\gamma |B(0,1)|}\Big)^{\frac{\ln(D_n(\eps, m, R) e^{n+2} 2^n |B(0,1)|^{\frac{1}{2}}(n+2)!)}{\ln 2}+\frac{\delta_{n,\eps,R}^{\eps-2}}{\ln 2} N^{1-\frac{\eps}{2}}}
\|f\|_{L^{\infty}(\omega \cap B_k)}.
\end{equation}
Recalling from the property (\ref{thick_rho}) that 
\begin{equation}{\label{thick2}}
|\omega \cap B_k| \geq \gamma |B_k|>0
\end{equation}
 as $B_k=B(x_k, \rho(x_k))$,  and
setting 
\begin{equation}\label{gh56}
\tilde{\omega}_{k}=\Big\{x \in \omega \cap B_k :\ |f(x)| \leq \frac{2}{|\omega \cap B_k|}\int_{\omega \cap B_k}|f(y)|dy\Big\},
\end{equation}
we observe that 
\begin{equation}\label{gh57}
\int_{\omega \cap B_k}|f(x)|dx \geq \int_{(\omega \cap B_k)\setminus \tilde{\omega}_{k}}|f(x)|dx \geq \frac{2|(\omega \cap B_k)\setminus \tilde{\omega}_{k}|}{|\omega \cap B_k|}\int_{\omega \cap B_k}|f(x)|dx.
\end{equation}
By using that the integral 
$$\int_{\omega \cap B_k}|f(x)|dx > 0,$$
is positive\footnote{This property can also be seen as a consequence of the Remez inequality, see e.g.~\cite{kkj} (Section~4.4).}, since $f$ is a non-zero entire function and $|\omega \cap B_k|>0$, 
we obtain that
$$|(\omega \cap B_k)\setminus \tilde{\omega}_{k}| \leq \frac{1}{2}|\omega \cap B_k|,$$
which implies that 
\begin{equation}\label{gh58}
|\tilde{\omega}_{k}|=|\omega \cap B_k|- |(\omega \cap B_k)\setminus \tilde{\omega}_{k}| \geq \frac{1}{2}|\omega \cap B_k|\geq \frac{1}{2} \gamma|B_k|= \frac{1}{2}\gamma \rho(x_k)^n \left|B(0,1)\right|>0,
\end{equation}
thanks to (\ref{thick2}).
By using again spherical coordinates as in (\ref{gh9}) and (\ref{gh10}), we observe that 
\begin{multline}\label{gh9a}
|\tilde{\omega}_{k}|=|\tilde{\omega}_{k} \cap B_k|
=\big( 2\rho(x_k) \big)^n \int_0^{1}\Big(\int_{\mathbb{S}^{n-1}}\un_{\tilde{\omega}_{k}  \cap B_k}(y_k+ 2\rho(x_k)r \sigma)d\sigma\Big)r^{n-1}dr  \\
\leq \big( 2\rho(x_k) \big)^n\int_{\mathbb{S}^{n-1}}|\tilde{I}_{\sigma}|d\sigma,
\end{multline}
where
\begin{equation}\label{gh11a}
\tilde{I}_{\sigma}=\{r \in [0,1] :\ y_k+ 2\rho(x_k)r \sigma \in \tilde{\omega}_{k} \cap B_k\}.
\end{equation}
As in (\ref{gh12}), the estimate (\ref{gh9a}) implies that there exists $\tilde{\sigma}_0(k) \in \mathbb{S}^{n-1}$ such that 
\begin{equation}\label{gh12a}
|\tilde{\omega}_{k}| \leq \big( 2\rho(x_k) \big)^n|\mathbb{S}^{n-1}| |\tilde{I}_{\tilde{\sigma}_0(k)}|.
\end{equation}
We deduce from (\ref{gh58}) and (\ref{gh12a}) that 
\begin{equation}\label{gh13a}
1 \geq|\tilde{I}_{\tilde{\sigma}_0(k)}| \geq \frac{\gamma |B(0,1)|}{2^{n+1}|\mathbb{S}^{n-1}|}>0.
\end{equation}
Applying anew Lemma~\ref{lem_kov} with $I=[0,1]$, $E=\tilde{I}_{\tilde{\sigma}_0(k)} \subset [0,1]$ verifying $|E|=|\tilde{I}_{\tilde{\sigma}_0(k)}|>0$, and the analytic function $\Phi=\phi$ defined in (\ref{gh13b}) with $\sigma_0(k)$ replaced by $\tilde{\sigma}_0(k)$ satisfying $|\phi(0)| \geq 1$, we obtain that
\begin{multline}\label{gh20a}
\left|B_k\right|^{\frac{1}{2}}\frac{\sup_{x \in [0,1]}|f(y_k+2\rho(x_k) x \tilde{\sigma}_0(k))|}{\|f\|_{L^2(B_k)}} \\
\leq \Big(\frac{C}{|\tilde{I}_{\tilde{\sigma}_0(k)}|}\Big)^{\frac{\ln M}{\ln 2}}\left|B_k\right|^{\frac{1}{2}}\frac{\sup_{x \in \tilde{I}_{\tilde{\sigma}_0(k)}}|f(y_k+2\rho(x_k) x \tilde{\sigma}_0(k))|}{\|f\|_{L^2(B_k)}},
\end{multline}
where $M \geq 1$ denotes the constant defined in (\ref{gh22}).
It follows from (\ref{gh13a}) and (\ref{gh20a}) that 
\begin{align}\label{gh21a}
& \ \sup_{x \in [0,1]}|f(y_k+2\rho(x_k) x \tilde{\sigma}_0(k))| \\ \notag
\leq & \ \Big(\frac{2^{n+1} C |\mathbb{S}^{n-1}|}{\gamma |B(0,1)|}\Big)^{\frac{\ln M}{\ln 2}}\sup_{x \in \tilde{I}_{\tilde{\sigma}_0(k)}}|f(y_k+2\rho(x_k) x \tilde{\sigma}_0(k))|\\ \notag
\leq & \ M^{\frac{1}{\ln 2}\ln(\frac{2^{n+1} C |\mathbb{S}^{n-1}|}{\gamma |B(0,1)|})}
\sup_{x \in \tilde{I}_{\tilde{\sigma}_0(k)}}|f(y_k+2\rho(x_k) x \tilde{\sigma}_0(k))|.
\end{align}
According to (\ref{gh11a}), we notice that 
\begin{equation}\label{gh23a}
\sup_{x \in \tilde{I}_{\tilde{\sigma}_0(k)}}|f(y_k+2\rho(x_k) x \tilde{\sigma}_0(k))| \leq \|f\|_{L^{\infty}(\tilde{\omega}_{k} \cap B_k)}.
\end{equation}
It follows from (\ref{gh8}), (\ref{gh21a}) and (\ref{gh23a}) that 
\begin{multline}\label{gh25a}
\|f\|_{L^{\infty}(B_k)}=|f(y_k)| \leq \sup_{x \in [0,1]}|f(y_k+2\rho(x_k) x \tilde{\sigma}_0(k))| \\
\leq M^{\frac{1}{\ln 2}\ln(\frac{2^{n+1} C |\mathbb{S}^{n-1}|}{\gamma |B(0,1)|})}\|f\|_{L^{\infty}(\tilde{\omega}_{k} \cap B_k)}.
\end{multline}
On the other hand, it follows from (\ref{gh56}) that
\begin{equation}\label{gh60}
\|f\|_{L^{\infty}(\tilde{\omega}_{k} \cap B_k)} \leq \frac{2}{|\omega \cap B_k|}\int_{\omega \cap B_k}|f(x)|dx.
\end{equation}
We deduce from (\ref{gh25a}), (\ref{gh60}) and the Cauchy-Schwarz inequality that 
\begin{align}\label{gh61}
 \|f\|_{L^{2}(B_k)} \leq & \ |B_k|^{\frac{1}{2}} \|f\|_{L^{\infty}(B_k)}\\ \notag
\leq & \ \frac{2|B_k|^{\frac{1}{2}}}{|\omega \cap B_k|} M^{\frac{1}{\ln 2}\ln(\frac{2^{n+1} C |\mathbb{S}^{n-1}|}{\gamma |B(0,1)|})}
\int_{\omega \cap B_k}|f(x)|dx\\ \notag
\leq & \  \frac{2|B_k|^{\frac{1}{2}}}{|\omega \cap B_k|^{\frac{1}{2}}} M^{\frac{1}{\ln 2}\ln(\frac{2^{n+1} C |\mathbb{S}^{n-1}|}{\gamma |B(0,1)|})}
\|f\|_{L^2(\omega \cap B_k)}.
\end{align}
By using the property (\ref{thick2}), it follows from (\ref{gh52}), (\ref{gh52s}) and (\ref{gh61}) that
\begin{multline}\label{gh62}
\|f\|_{L^{2}(B_k)}^2 \leq  \frac{4}{\gamma} M^{\frac{2}{\ln 2}\ln(\frac{2^{n+1} C |\mathbb{S}^{n-1}|}{\gamma |B(0,1)|})}
\|f\|_{L^2(\omega \cap B_k)}^2\\
\leq  \frac{4}{\gamma} \Big(|B(0,1)|^{\frac{1}{2}} (n+2)! D_n(\eps, m, R) e^{n+2} 2^n e^{\delta_{n,\eps,R}^{\eps-2} N^{1-\frac{\eps}{2}}}\Big)^{\frac{2}{\ln 2}\ln(\frac{2^{n+1} C |\mathbb{S}^{n-1}|}{\gamma |B(0,1)|})}
\|f\|_{L^2(\omega \cap B_k)}^2.
\end{multline}
Setting
\begin{equation}\label{gh54}
\kappa_n(m, R,\gamma,\eps)=\frac{4}{\sqrt{3}\gamma^{\frac{1}{2}}}
\Big(\frac{2^{n+1} C |\mathbb{S}^{n-1}|}{\gamma |B(0,1)|}\Big)^{\frac{\ln(|B(0,1)|^{\frac{1}{2}} (n+2)! D_n(\eps, m, R) e^{n+2} 2^n)}{\ln 2}}>0,
\end{equation}
we deduce from (\ref{gh62}) that there exists a positive universal constant $\tilde{\kappa}_n>1$ such that for any good ball $B_k$,
\begin{equation}\label{gh55}
\|f\|_{L^{2}(B_k)}^2\leq \frac{3}{4}\kappa_n(m, R,\gamma,\eps)^2\Big(\frac{\tilde{\kappa}_n}{\gamma}\Big)^{\frac{2}{\ln 2}\delta_{n,\eps,R}^{\eps-2} N^{1-\frac{\eps}{2}}}
\|f\|_{L^{2}(\omega \cap B_k)}^2.
\end{equation}
By using anew from (\ref{recov}) that 
\begin{equation}
\mathbbm{1}_{\bigcup_{\textrm{good balls}} B_k}  \leq \sum \limits_{\textrm{good balls}} \mathbbm{1}_{B_k} \leq N_0 \mathbbm{1}_{\bigcup_{\textrm{good balls}} B_k},
\end{equation}
it follows from (\ref{gh7}) and (\ref{gh55}) that 
\begin{align}\label{gh56y}
\|f\|_{L^2(\rr^n)}^2 & \leq \frac{4}{3} \int_{\bigcup_{\textrm{good balls}} B_k}|f(x)|^2dx \leq \frac{4}{3} \sum_{\textrm{good balls}}\|f\|_{L^{2}(B_k)}^2\\
\notag
& \leq \kappa_n(m, R,\gamma, \eps)^2\Big(\frac{\tilde{\kappa}_n}{\gamma}\Big)^{\frac{2}{\ln 2}\delta_{n,\eps,R}^{\eps-2} N^{1-\frac{\eps}{2}}} \sum_{\textrm{good balls}}
\|f\|_{L^{2}(\omega \cap B_k)}^2\\
\notag
& \leq N_0 \kappa_n(m, R,\gamma, \eps)^2\Big(\frac{\tilde{\kappa}_n}{\gamma}\Big)^{\frac{2}{\ln 2}\delta_{n,\eps,R}^{\eps-2} N^{1-\frac{\eps}{2}}} \int_{\omega \cap (\bigcup_{\textrm{good balls}} B_k)}|f(x)|^2dx\\ 
\notag
& \leq N_0 \kappa_n(m, R,\gamma,\eps)^2\Big(\frac{\tilde{\kappa}_n}{\gamma}\Big)^{\frac{2}{\ln 2}\delta_{n,\eps,R}^{\eps-2}N^{1-\frac{\eps}{2}}}\|f\|_{L^2(\omega)}^2.
\end{align}
This ends the proof of Theorem~\ref{Spectral}.

\section{Proof of Theorem~\ref{control_PDEfrac}}\label{control_PDEfrac_proof}
This section is devoted to the proof of Theorem~\ref{control_PDEfrac}. Let $A$ be a closed operator on $L^2(\rr^n)$ which is the  infinitesimal generator of a strongly continuous contraction semigroup $(e^{-tA})_{t \geq 0}$ on $L^2(\rr^n)$ satisfying the assumptions of Theorem~\ref{control_PDEfrac}. According to Lemma~\ref{equivv}, we can assume that there exist some constants 
$\frac{1}{2}<s \leq 1$, $C_s > 1$, $0<t_0 \leq 1$, $m_1,m_2 \in \rr$ with $m_1>0$, $m_2 \geq 0$  such that
\begin{multline}\label{GS0_bis}
\forall 0< t \leq t_0, \forall \alpha, \beta  \in \nn^n, \forall g \in L^2(\rr^n),\\
\| x^{\alpha} \partial_x^{\beta}( e^{-t A}g)\|_{L^2(\rr^n)} \leq  \frac{C_s^{1+|\alpha|+|\beta|}}{t^{m_1 (|\alpha|+|\beta|)+m_2}} (\alpha !)^{\frac{1}{2s}}(\beta!)^{\frac{1}{2s}}  \|g\|_{L^2(\rr^n)}.
\end{multline}
Let $\rho : \rr^n \longrightarrow (0,+\infty)$ be a $\frac{1}{2}$-Lipschitz positive function
with~$\rr^n$ being equipped with the Euclidean norm,
satisfying that there exist some constants $0 \leq  \delta < 2s-1$, $m>0$, $R>0$ such that
\begin{equation*}
\forall x \in \rr^n, \quad 0<m \leq \rho(x) \leq R{\left\langle x\right\rangle}^{\delta}.
\end{equation*} 
Let $\omega$ be a measurable subset of $\rr^n$ which is $\gamma$-thick with respect to the density $\rho$, that is, 
\begin{equation}{\label{thick_rho0_bis}}
 \exists 0< \gamma \leq 1, \forall x \in \rr^n, \quad |\omega \cap B(x,\rho(x))| \geq \gamma |B(x,\rho(x))|.
\end{equation}
Thanks to the Hilbert Uniqueness Method, the null-controllability of the system (\ref{ucla1}) is equivalent to the observability of the adjoint system 
\begin{equation*}
\left\lbrace \begin{array}{ll}
(\partial_t +A)g(t,x)=0\,, \quad & x \in \mathbb{R}^n,\ t>0, \\
g|_{t=0}=g_0 \in L^2(\rr^n),
\end{array}\right.
\end{equation*}
from the control subset $\omega$ in any positive time $T>0$. We shall prove that Theorem~\ref{control_PDEfrac} can be deduced from the abstract observability result given by Theorem~\ref{Meta_thm_AdaptedLRmethod}. In order to apply Theorem~\ref{Meta_thm_AdaptedLRmethod}, it is therefore sufficient to check that the spectral inequality (\ref{Meta_thm_IS}) and the dissipation estimate (\ref{Meta_thm_dissip}) hold when using the Hermite orthogonal projections $(\pi_k)_{k \geq 0}$ defined in (\ref{projection}).
It follows from Theorem~\ref{Spectral} that 
there exist some positive constant $\kappa_n(m, R, \gamma,1-\delta)>0$, $\tilde{C}_n(1-\delta, R) >0$ and a positive universal constant $\tilde{\kappa}_n >0$ such that for all $k \geq 1$, $f \in L^2(\rr^n)$,  
\begin{equation}\label{ucla2}
\|\pi_kf\|_{L^2(\rr^n)} \leq \kappa_n(m, R, \gamma,1-\delta) \Big( \frac{\tilde{\kappa}_n}{\gamma} \Big)^{\tilde{C}_n(1-\delta, R) k^{\frac{1+\delta}{2}}} \|\pi_kf\|_{L^2(\omega)}.
\end{equation}
This establishes the spectral inequality (\ref{Meta_thm_IS}) with the parameter $0<a=\frac{1+\delta}{2}<s$. Let us now prove that the dissipation estimate (\ref{Meta_thm_dissip}) holds true as well. To that end, we begin by establishing that there exists a positive constant $\tilde{C}_s(n) >1$ such that for all $k \in \nn$, $g \in L^2(\rr^n)$, $ 0 < t \leq t_0$,
\begin{equation}\label{GS02}
\|(\mathcal{H}+n)^{k}(e^{-tA}g) \|_{L^2(\rr^n)} \leq  \frac{\tilde{C}_s(n)^{1+k}}{t^{2m_1 k +m_2}} (k!)^{\frac{1}{s}} \|g\|_{L^2(\rr^n)},
\end{equation}
where $\mathcal{H}= \sum_{j=1}^n \mathcal{H}_j$ denotes the harmonic oscillator with 
\begin{equation}\label{deft0}
\mathcal{H}_j+1= -\partial^2_{x_j} + x_j^2+1=(\partial_{x_j}+x_j)(-\partial_{x_j} +x_j), \qquad 1\leq j \leq n.
\end{equation} 
Let $k \in \nn^*$. We deduce from (\ref{deft0}) and Lemma~\ref{dvpt} that there exists a finite family of real numbers $(C_{l_1,l_2}^{2k-1})_{\substack{l_1, l_2 \in \nn, \\ 0 \leq l_1+ l_2 \leq 2k}}$ independent on the parameter $1 \leq j \leq n$ such that
\begin{equation}\label{deft1}
(\mathcal{H}_j +1)^k=\sum_{\substack{l_1,l_2 \in \nn, \\ 0 \leq l_1 + l_2 \leq 2k,}} C^{2k-1}_{l_1,l_2} x_j^{l_1} \partial_{x_j}^{l_2}.
\end{equation}
and
\begin{equation}\label{deft2}
\forall l_1, l_2 \in \nn, \ 0 \leq l_1+ l_2 \leq 2k, \quad |C^{2k-1}_{l_1, l_2}| \leq 3^{2k-1} (2k)^{\frac{2k-l_1-l_2}{2}}.
\end{equation}
By using that $[\mathcal{H}_j+1,\mathcal{H}_k+1]=0$ for all $0 \leq j,k \leq n$,
we deduce from the multinomial formula that 
\begin{align}\label{har1}
& \ (\mathcal{H}+n)^k\\ \notag
= & \ \sum_{\substack{\gamma=(\gamma_1,...,\gamma_n) \in \nn^n, \\ |\gamma|=k}} \frac{k!}{\gamma!} \prod_{j=1}^n (\mathcal{H}_j+1)^{\gamma_j}
=\sum_{\substack{\gamma=(\gamma_1,...,\gamma_n) \in \nn^n, \\ |\gamma|=k}} \frac{k!}{\gamma!} 
\prod_{j=1}^n \sum_{\substack{l_j,\tilde{l}_j \in \nn, \\ 0 \leq l_j + \tilde{l}_j \leq 2\gamma_j}} C^{2\gamma_j-1}_{l_j, \tilde{l}_j} x_j^{l_j} \partial_{x_j}^{\tilde{l}_j} \\ \notag
= & \ \sum_{\substack{\gamma \in \nn^n, \\ |\gamma|=k}} \frac{k!}{\gamma!} \sum_{\substack{\alpha, \beta \in \nn^n, \\ \alpha + \beta \leq 2\gamma}} \big(\prod_{j=1}^n C^{2\gamma_j-1}_{\alpha_j, \beta_j} \big) x^{\alpha} \partial_x^{\beta} = \sum_{\substack{\alpha, \beta \in \nn^n, \\ |\alpha + \beta| \leq 2k}} \sum_{\substack{\gamma \in \nn^n,  |\gamma|=k, \\  \alpha + \beta \leq 2\gamma}} \frac{k!}{\gamma!} \big(\prod_{j=1}^n C^{2\gamma_j-1}_{\alpha_j, \beta_j} \big) x^{\alpha} \partial_x^{\beta}.
\end{align}
It follows from \eqref{har1} that
\begin{equation}\label{deft4}
(\mathcal{H}+n)^k = \sum_{\substack{\alpha, \beta \in \nn^n, \\ |\alpha+\beta| \leq 2k}} c^k_{\alpha, \beta} x^{\alpha} \partial_x^{\beta},
\end{equation}
with 
\begin{equation}\label{deft5}
c^k_{\alpha, \beta} = \sum_{\substack{\gamma \in \nn^n, |\gamma|=k, \\  \alpha + \beta \leq 2\gamma}} \frac{k!}{\gamma!} \big(\prod_{j=1}^n C^{2\gamma_j-1}_{\alpha_j, \beta_j} \big).
\end{equation}
It follows from (\ref{deft2}) and (\ref{deft5}) that for all $\alpha, \beta \in \nn^n$ with $|\alpha+\beta| \leq 2k$,
\begin{multline}\label{har2}
|c^k_{\alpha, \beta}| \leq \sum_{\substack{\gamma \in \nn^n, |\gamma|=k, \\  \alpha + \beta \leq 2\gamma}} \frac{k!}{\gamma!} \Big(\prod_{j=1}^n 3^{2\gamma_j-1} (2\gamma_j)^{\frac{2\gamma_j-\alpha_j-\beta_j}{2}}\Big)\\
 \leq \sum_{\substack{\gamma \in \nn^n, \\ |\gamma|=k}} \frac{k!}{\gamma!} 3^{2|\gamma|-n} (2k)^{\frac{2k-|\alpha+\beta|}{2}} =  3^{2k-n}n^k (2k)^{\frac{2k-|\alpha+\beta|}{2}}.
\end{multline}
We deduce from (\ref{GS0_bis}), (\ref{deft4}) and \eqref{har2} that for all $k \geq 1$, $g \in L^2(\rr^n)$, $0< t \leq t_0$, 
\begin{align}\label{deft8}
& \ \|(\mathcal{H}+n)^{k}(e^{-tA}g)\|_{L^2(\rr^n)} \leq \sum_{\substack{\alpha, \beta \in \nn^n, \\ |\alpha + \beta| \leq 2k}} |c^{k}_{\alpha, \beta}| \|x^{\alpha} \partial^{\beta}_x(e^{-tA}g)\|_{L^2(\rr^n)} \\  \notag
\leq & \ \sum_{\substack{\alpha, \beta \in \nn^n, \\ |\alpha + \beta| \leq 2k}} 3^{2k-n}n^k (2k)^{\frac{2k-|\alpha| -|\beta|}{2}}  \frac{C_s^{1+|\alpha| + |\beta|}}{t^{m_1(|\alpha|+ |\beta|)+ m_2}} (\alpha!)^{\frac{1}{2s}}(\beta!)^{\frac{1}{2s}} \|g\|_{L^2(\rr^n)} \\ \notag
\leq & \ \sum_{\substack{\alpha, \beta \in \nn^n, \\ |\alpha + \beta| \leq 2k}}  3^{2k-n}n^k (2k)^{s \frac{2k-|\alpha|-|\beta|}{2s}}  \frac{C_s^{1+|\alpha| + |\beta|}}{t^{m_1(|\alpha|+ |\beta|)+ m_2}} (|\alpha|)^{\frac{|\alpha|}{2s}}(|\beta|)^{\frac{|\beta|}{2s}} \|g\|_{L^2(\rr^n)},
\end{align}
when using above the convention $0^0=1$. With this convention, we directly notice from Lemma~\ref{gamma3} that 
$$\forall x, y \geq 0, \quad x^y \leq e^x \Gamma(y+1) \quad \text{and} \quad \forall x, y \geq 1, \quad \Gamma(x) \Gamma(y) \leq \frac{B(1,1)}{2} \Gamma(x+y+1).$$
By using that $\frac{1}{2}<s \leq 1$ and the above estimates, it follows from (\ref{deft8}) that for all $k \geq 1$, $g \in L^2(\rr^n)$, $0< t \leq t_0$,
\begin{align}\label{GS3}
& \ \|(\mathcal{H}+n)^{k}(e^{-tA}g)\|_{L^2(\rr^n)} \\ \notag
\leq & \  \sum_{\substack{\alpha, \beta \in \nn^n, \\ |\alpha + \beta| \leq 2k}} 3^{2k-n}n^k e^{(2k)^s} \Gamma\Big(\frac{2k-|\alpha|-|\beta|}{2s}+1 \Big)  \frac{C_s^{1+|\alpha| + |\beta|}e^{|\alpha+ \beta|}}{t^{m_1(|\alpha|+ |\beta|)+ m_2}}  \Gamma\Big(\frac{|\alpha|}{2s}+1\Big) \Gamma\Big(\frac{|\beta|}{2s}+1\Big) \|g\|_{L^2} \\ \notag
\leq & \ \sum_{\substack{\alpha, \beta \in \nn^n, \\ |\alpha + \beta| \leq 2k}} 3^{2k-n}n^k e^{4k} \Gamma\Big(\frac{2k-|\alpha|-|\beta|}{2s}+1 \Big)  \frac{C_s^{1+|\alpha| + |\beta|}}{t^{m_1(|\alpha|+ |\beta|)+ m_2}} \Gamma\Big(\frac{|\alpha|}{2s}+1\Big) \Gamma\Big(\frac{|\beta|}{2s}+1\Big) \|g\|_{L^2} \\ \notag
\leq &\  \sum_{\substack{\alpha, \beta \in \nn^n,\\ |\alpha + \beta| \leq 2k}} 3^{2k-n} n^k e^{4k}\frac{B(1,1)^2}{4} \Gamma\Big(\frac{k}{s}+5 \Big)  \frac{C_s^{1+|\alpha| + |\beta|}}{t^{m_1(|\alpha|+ |\beta|)+ m_2}}   \|g\|_{L^2}.
\end{align}
Thanks to Stirling formula (\ref{stirling1}), we can find a positive constant $C'_s >1$ such that for all $k \geq 1$, 
\begin{equation}\label{stirling}
\Gamma\Big(\frac{k}{s}\Big) \leq C'_s \sqrt{\frac{2\pi s}{k}} s^{-\frac{k}{s}} e^{-\frac{k}{s}} k^{\frac{k}{s}} \leq C'_s \sqrt{\frac{2\pi s}{k}} s^{-\frac{k}{s}} (k!)^{\frac{1}{s}},
\end{equation}
since 
$$\forall k \geq 1, \quad \frac{k^k}{k!} \leq \sum_{j=0}^{+\infty}\frac{k^j}{j!}=e^k.$$
By using from (\ref{gamma2}) that 
$$\Gamma\Big(\frac{k}{s}+5\Big)= \Big(\frac{k}{s}+4\Big)\Big(\frac{k}{s}+3\Big)\Big(\frac{k}{s}+2\Big)\Big(\frac{k}{s}+1\Big)\frac{k}{s}\Gamma\Big(\frac{k}{s}\Big),$$ 
it follows from \eqref{GS3} and \eqref{stirling} that for all $k \geq 1$, $g \in L^2(\rr^n)$, $0 < t \leq t_0$,
\begin{align}\label{deft11}
&\ \|(\mathcal{H}+n)^{k}(e^{-tA}g)\|_{L^2(\rr^n)} \\ \notag
\leq & \ \sum_{\substack{\alpha,\beta \in \nn^n, \\ |\alpha + \beta| \leq 2k}} 3^{2k-n} n^k e^{4k} \frac{B(1,1)^2}{4}\Gamma\Big(\frac{k}{s}+5 \Big)   \frac{C_s^{1+|\alpha| + |\beta|}}{t^{m_1(|\alpha|+ |\beta|)+ m_2}}   \|g\|_{L^2(\rr^n)} \\ \notag
\leq &\ \sum_{\substack{\alpha,\beta \in \nn^n, \\ |\alpha + \beta| \leq 2k}}  3^{2k-n} n^ke^{4k} \frac{B(1,1)^2}{4} \Big(\frac{k}{s}+4\Big)^5 C'_s \sqrt{\frac{2\pi s}{k}} s^{-\frac{k}{s}}   \frac{C_s^{1+|\alpha| + |\beta|}}{t^{m_1(|\alpha|+ |\beta|)+ m_2}} (k!)^{\frac{1}{s}} \|g\|_{L^2(\rr^n)} \\ \notag
\leq &\  (2k+1)^{2n}  3^{2k-n} n^k e^{4k} \frac{B(1,1)^2}{4} \Big(\frac{k}{s}+4\Big)^5 C'_s \sqrt{\frac{2\pi s}{k}} s^{-\frac{k}{s}}   \frac{C_s^{1+2k}}{t^{2m_1 k+ m_2}}  (k!)^{\frac{1}{s}} \|g\|_{L^2(\rr^n)}. 
\end{align}
We deduce from (\ref{deft11}) that there exists a positive constant $\tilde{C}_s(n)>1$ such that for all $k \geq 1$, $g \in L^2(\rr^n)$, $0 < t \leq t_0$,
\begin{equation}\label{esti2}
\|(\mathcal{H}+n)^{k}(e^{-tA}g)\|_{L^2(\rr^n)} \leq    \frac{\tilde{C}_s(n)^{1+k}}{t^{2m_1 k+m_2}} (k!)^{\frac{1}{s}} \|g\|_{L^2(\rr^n)}.
\end{equation}
The estimate (\ref{esti2}) holds as well when $k=0$ since $(e^{-tA})_{t \geq 0}$ is a contraction semigroup on $L^2(\rr^n)$, $\tilde{C}_s(n)>1$ and $0<t_0 \leq 1$. This ends the proof of the estimates (\ref{GS02}).

With $(\Phi_{\alpha})_{\alpha \in \nn^n}$ the $L^2(\rr^n)$-Hermite basis, we next notice that for all $g \in L^2(\rr^n)$, $t \geq 0$,
\begin{align}\label{de1}
& \ \sum_{\alpha \in \nn^n}e^{2t^{2m_1s+1}(2|\alpha|+n)^s} |\langle e^{-tA}g, \Phi_{\alpha}\rangle_{L^2(\rr^n)} |^2  \\ \notag
= & \ \sum_{\alpha \in \nn^n} \sum_{k=0}^{+\infty} \frac{2^k t^{k(2m_1s+1)} (2|\alpha|+ n)^{sk}}{k!}|\langle e^{-tA}g, \Phi_{\alpha}\rangle_{L^2(\rr^n)}|^2  \\ \notag
 \leq & \ \sum_{\alpha \in \nn^n} \sum_{k=0}^{+\infty} \frac{2^k t^{k(2m_1s+1)}}{k!} |\langle e^{-tA}g, (2|\alpha|+2n)^{\lfloor \frac{ks}{2} \rfloor +1} \Phi_{\alpha}\rangle_{L^2(\rr^n)} |^2, 
\end{align}
where $\lfloor \cdot \rfloor$ denotes the floor function. 
By using the selfadjointness property of the harmonic oscillator $\mathcal{H}=-\Delta_x+|x|^2$, we deduce from (\ref{esti2}), (\ref{de1}) and (\ref{6.harmo})  that for all $g \in L^2(\rr^n)$, $t \geq 0$,
\begin{align}\label{de2}
& \ \sum_{\alpha \in \nn^n}e^{2t^{2m_1s+1}(2|\alpha|+n)^s} |\langle e^{-tA}g, \Phi_{\alpha}\rangle_{L^2(\rr^n)} |^2 \\ \notag
 \leq & \  \sum_{\alpha \in \nn^n} \sum_{k=0}^{+\infty} \frac{2^k t^{k(2m_1s+1)}}{k!} |\langle e^{-tA}g,(\mathcal{H}+n)^{\lfloor \frac{ks}{2} \rfloor +1} \Phi_{\alpha}\rangle_{L^2(\rr^n)}|^2 \\ \notag
= &\ \sum_{k=0}^{+\infty} \frac{2^k t^{k(2m_1s+1)}}{k!}  \sum_{\alpha \in \nn^n} |\langle (\mathcal{H}+n)^{\lfloor \frac{ks}{2} \rfloor +1}(e^{-tA}g), \Phi_{\alpha}\rangle_{L^2(\rr^n)}|^2\\ \notag 
= & \ \sum_{k=0}^{+\infty} \frac{2^k t^{k(2m_1s+1)}}{k!} \|(\mathcal{H}+n)^{\lfloor \frac{ks}{2} \rfloor +1}(e^{-tA}g) \|^2_{L^2(\rr^n)}.
\end{align}
It follows from (\ref{GS02}) and (\ref{de2}) that for all $g \in L^2(\rr^n)$, $0<t \leq t_0$,
\begin{align}\label{deft20}
& \ \sum_{\alpha \in \nn^n}e^{2t^{2m_1s+1}(2|\alpha|+n)^s} |\langle e^{-tA}g, \Phi_{\alpha}\rangle_{L^2(\rr^n)} |^2  \\ \notag
& \leq \sum_{k=0}^{+\infty} \frac{2^k t^{k(2m_1s+1)}}{k!} \frac{\tilde{C}_s(n)^{2 \lfloor \frac{ks}{2} \rfloor +4}}{t^{4m_1 (\lfloor \frac{ks}{2} \rfloor +1)+2m_2}} \Big(\Big(\Big\lfloor \frac{ks}{2} \Big\rfloor +1\Big)!\Big)^{\frac{2}{s}}  \|g \|^2_{L^2(\rr^n)}.
\end{align}
By using that $\lfloor \frac{ks}{2} \rfloor \leq \frac{ks}{2}$ and that the Gamma function is increasing on $[2, +\infty)$, see Section~\ref{miscgamma}, we deduce from (\ref{deft20}) and (\ref{factoriel}) that for all $g \in L^2(\rr^n)$, $0< t \leq t_0$,
\begin{align}\label{GS2}
& \ \sum_{\alpha \in \nn^n}e^{2t^{2m_1s+1}(2|\alpha|+n)^s} |\langle e^{-tA}g, \Phi_{\alpha}\rangle_{L^2(\rr^n)} |^2  \\ \notag
& \leq \sum_{k=0}^{+\infty} \frac{2^k t^{k(2m_1s+1)}}{k!} \frac{\tilde{C}_s(n)^{4+ks}}{t^{2m_1 ks +4m_1+2m_2}} \Big(\Big(\Big\lfloor \frac{ks}{2} \Big\rfloor +1\Big)!\Big)^{\frac{2}{s}}  \|g \|^2_{L^2(\rr^n)} \\ \notag
& \leq \sum_{k=0}^{+\infty} \frac{2^k t^{k}}{k!} \frac{\tilde{C}_s(n)^{4+ks}}{t^{4m_1 +2m_2}} \Gamma\Big(\Big\lfloor \frac{ks}{2} \Big\rfloor +2\Big)^{\frac{2}{s}}  \|g \|^2_{L^2(\rr^n)} \\ \notag
& \leq \sum_{k=0}^{+\infty} \frac{2^k t^{k}}{k!} \frac{\tilde{C}_s(n)^{4+ks}}{t^{4m_1 +2m_2}} \Gamma\Big(\frac{ks}{2} +2\Big)^{\frac{2}{s}}  \|g \|^2_{L^2(\rr^n)}.
\end{align}
By using Lemma~\ref{gamma3} (assertion $(iii)$) and $\frac{2}{s} \geq 1$ since $\frac{1}{2}<s \leq 1$, we can find a positive constant $C''_s>1$ such that 
\begin{equation}\label{deft30}
\forall x\geq 1, \quad  \Gamma(x)^{\frac{2}{s}} \leq C''_s \Gamma\Big(\frac{2x}{s}\Big). 
\end{equation}
With the notation $\tilde{C}_s=\tilde{C}_s(n)$, we deduce from (\ref{GS2}), (\ref{deft30}) and (\ref{factoriel}) that for all $g \in L^2(\rr^n)$, $0< t \leq t_0$,
\begin{align}\label{deft31}
& \ \sum_{\alpha \in \nn^n}e^{2t^{2m_1s+1}(2|\alpha|+n)^s} |\langle e^{-tA}g, \Phi_{\alpha}\rangle_{L^2} |^2 
 \leq  C''_s \sum_{k=0}^{+\infty} \frac{2^k t^{k}}{k!} \frac{\tilde{C}_s^{4+ks}}{t^{4m_1 +2m_2}} \Gamma\Big(k +\frac{4}{s}\Big)  \|g \|^2_{L^2} \\ \notag
\leq &  \ C''_s \sum_{k=0}^{+\infty} \frac{2^k t^{k}}{k!} \frac{\tilde{C}_s^{4+ks}}{t^{4m_1 +2m_2}} \Gamma\big(k +8\big)  \|g \|^2_{L^2} 
 \leq  C''_s \sum_{k=0}^{+\infty} \frac{2^k t^{k}}{k!} \frac{\tilde{C}_s^{4+ks}}{t^{4m_1 +2m_2}} (k +7)!  \|g \|^2_{L^2} \\ \notag
 \leq & \ C''_s \sum_{k=0}^{+\infty} 2^k t^{k} \frac{\tilde{C}_s^{4+ks}}{t^{4m_1 +2m_2}} (k +7)^7  \|g \|^2_{L^2} 
\leq  C''_s \Big(\sum_{k=0}^{+\infty}  (2\tilde{C}_s^s e t)^k\Big)\frac{\tilde{C}_s^{4}}{t^{4m_1 +2m_2}} 7! e^7  \|g \|^2_{L^2},
\end{align}
since $\frac{1}{2}<s \leq 1$ and
$$\frac{(k+7)^7}{7!} \leq \sum_{j=0}^{+\infty}\frac{(k+7)^j}{j!} \leq e^{k+7}.$$
It follows from \eqref{deft31} that for all $g \in L^2(\rr^n)$, $0< t \leq t_1$,
\begin{equation}\label{Hermite_decay3}
\sum_{\alpha \in \nn^n}e^{2t^{2m_1s+1}(2|\alpha|+n)^s} |\langle e^{-tA}g, \Phi_{\alpha}\rangle_{L^2(\rr^n)} |^2 \leq \frac{2 C''_s\tilde{C}_s^{4}}{t^{4m_1 +2m_2}} 7! e^7  \|g \|^2_{L^2(\rr^n)},
\end{equation}
with 
$$0<t_1=\min(t_0, (4\tilde{C}_s^s e)^{-1}) \leq 1.$$
For any $0 < t \leq t_1$ and $g \in L^2(\rr^n)$, the series 
\begin{equation*}
 f= \sum_{\alpha \in \nn^n} e^{t^{2m_1s+1}(2|\alpha|+n)^s}\langle e^{-tA}g, \Phi_\alpha \rangle_{L^2(\rr^n)} \Phi_\alpha,
 \end{equation*}
is therefore convergent in $L^2(\rr^n)$ and defines a $L^2(\rr^n)$-function satisfying
\begin{equation}\label{ucla20}
\|f\|_{L^2(\rr^n)} \leq \frac{\sqrt{2C_s'' 7!}\tilde{C}_s^2e^{\frac{7}{2}}}{t^{2m_1+m_2}} \|g\|_{L^2(\rr^n)}, \qquad e^{-t^{2m_1s+1}\mathcal{H}^s}f=e^{-tA}g,
 \end{equation}
 according to (\ref{SG_frac}).
It follows from (\ref{ucla20}) that for all $0 < t \leq t_1$, $g \in L^2(\rr^n)$, $k \geq 1$,
\begin{multline}{\label{Esti_dissip1}}
\|(1-\pi_k)(e^{-tA}g)\|_{L^2} = \|(1-\pi_k)(e^{-t^{2m_1s+1}\mathcal{H}^s}f)\|_{L^2} 
= \|e^{-t^{2m_1s+1}\mathcal{H}^s}(1-\pi_k)f\|_{L^2}\\
 \leq e^{-t^{2m_1s+1}(2k+n)^s} \|(1-\pi_k)f\|_{L^2}\leq e^{-t^{2m_1s+1}k^s} \|f\|_{L^2}. 
\end{multline}
We deduce from (\ref{ucla20}) and (\ref{Esti_dissip1}) the following dissipation estimate
\begin{multline}{\label{Esti_dissip2}}
\forall 0 < t \leq t_1, \forall g \in L^2(\rr^n), \forall k \geq 1,\\
\|(1-\pi_k)(e^{-tA}g)\|_{L^2(\rr^n)}  \leq \frac{\sqrt{2C_s'' 7!}\tilde{C}_s^2e^{\frac{7}{2}}}{t^{2m_1+m_2}} e^{-t^{2m_1s+1}k^s} \|g\|_{L^2(\rr^n)}.
\end{multline} 
It establishes the dissipation estimate (\ref{Meta_thm_dissip}) with the parameter $0<a=\frac{1+\delta}{2}<b=s$. We can therefore deduce from Theorem~\ref{Meta_thm_AdaptedLRmethod} that the following observability estimate holds in any positive time 
\begin{multline*}
\exists C>1, \forall T>0, \forall g \in L^2(\rr^n), \\ 
\|e^{-TA}g\|_{L^2(\rr^n)}^2 \leq C\exp\Big(\frac{C}{T^{\frac{(1+\delta)(2m_1s+1)}{2s-1-\delta}}}\Big) \int_0^T \|e^{-tA}g\|_{L^2(\omega)}^2 dt.
\end{multline*}
This ends the proof of Theorem~\ref{control_PDEfrac}.

We close this section by noticing that the conclusions of Theorem~\ref{control_PDEfrac} hold true as well when the quantitative regularizing estimates (\ref{GS0}) holding for some $\frac{1}{2}< s \leq1$ are replaced by the following assumption 
\begin{multline}{\label{decay_hermite_bis}}
\exists \frac{1}{2}<s \leq 1,\exists m_1, m_2>0, \exists C_1, C_2 >0, \exists 0<t_0\leq 1, \forall 0 < t \leq t_0, \forall g \in L^2(\rr^n),\\ \sum_{\alpha \in \nn^n}e^{\frac{2 t^{m_1}}{C_1} (2 |\alpha|+n)^s} |\langle e^{-t A}g, \Phi_\alpha\rangle_{L^2(\rr^n)}|^2  \leq \frac{C_2^2}{t^{2m_2}} \|g\|_{L^2(\rr^n)}^2. 
\end{multline}
By resuming the above proof from (\ref{Hermite_decay3}), we indeed notice that for any $0 < t \leq t_0$ and $g \in L^2(\rr^n)$, the series 
\begin{equation*}
 f= \sum_{\alpha \in \nn^n} e^{\frac{t^{m_1}}{C_1}(2|\alpha|+n)^s}\langle e^{-tA}g, \Phi_\alpha \rangle_{L^2(\rr^n)} \Phi_\alpha,
 \end{equation*}
is convergent in $L^2(\rr^n)$ and defines a $L^2(\rr^n)$-function satisfying
\begin{equation}\label{ucla20_bis}
\|f\|_{L^2(\rr^n)} \leq \frac{C_2}{t^{m_2}} \|g\|_{L^2(\rr^n)}, \qquad e^{-\frac{t^{m_1}}{C_1}\mathcal{H}^s}f=e^{-tA}g,
 \end{equation}
according to (\ref{SG_frac}).
It follows from (\ref{ucla20_bis}) that for all $0 < t \leq t_0$, $g \in L^2(\rr^n)$, $k \geq 1$,
\begin{multline}{\label{Esti_dissip1_bis}}
\|(1-\pi_k)(e^{-tA}g)\|_{L^2(\rr^n)} = \|(1-\pi_k)(e^{-\frac{t^{m_1}}{C_1}\mathcal{H}^s}f)\|_{L^2(\rr^n)} \\
= \|e^{-\frac{t^{m_1}}{C_1}\mathcal{H}^s}(1-\pi_k)f\|_{L^2(\rr^n)}
 \leq e^{-\frac{t^{m_1}}{C_1}(2k+n)^s} \|(1-\pi_k)f\|_{L^2(\rr^n)}\leq e^{-\frac{t^{m_1}}{C_1}k^s} \|f\|_{L^2(\rr^n)}. 
\end{multline}
We deduce from (\ref{ucla20_bis}) and (\ref{Esti_dissip1_bis}) the following dissipation estimate
\begin{multline}{\label{Esti_dissip2_bis}}
\forall 0 < t \leq t_0, \forall g \in L^2(\rr^n), \forall k \geq 1,\quad
\|(1-\pi_k)(e^{-tA}g)\|_{L^2}  \leq \frac{C_2}{t^{m_2}} e^{-\frac{t^{m_1}}{C_1}k^s} \|g\|_{L^2}.
\end{multline} 
It establishes the dissipation estimate (\ref{Meta_thm_dissip}) with the parameter $0<a=\frac{1+\delta}{2}<b=s$. We can therefore deduce from Theorem~\ref{Meta_thm_AdaptedLRmethod} that the following observability estimate holds in any positive time 
$$\exists C>1, \forall T>0, \forall g \in L^2(\rr^n), \quad
\|e^{-TA}g\|_{L^2(\rr^n)}^2 \leq C\exp\Big(\frac{C}{T^{\frac{(1+\delta)m_1}{2s-1-\delta}}}\Big) \int_0^T \|e^{-tA}g\|_{L^2(\omega)}^2 dt.$$

\section{Appendix}\label{appendix}

\subsection{Miscellaneous facts about the Gamma function}\label{miscgamma}

Let $\nn$ be the set of non-negative integers and $\mathbb{Z}_-$ be the set of non-positive integers. 
The Gamma function defined as
\begin{equation}{\label{gamma1}}
\forall x >0, \quad \Gamma(x)=\int_0^{+\infty}t^{x-1} e^{-t} dt>0,
\end{equation}
admits an unique analytic extension on $\cc \setminus{\Z_-}$ satisfying the functional identity
\begin{equation}{\label{gamma2}}
\forall z \in \cc \setminus{\Z_-},\quad  \Gamma(z+1)= z \Gamma(z),
\end{equation}
and interpolating the factorial function
\begin{equation}{\label{factoriel}}
\forall n \in \nn, \quad \Gamma(n+1)=n!.
\end{equation}
It also satisfies the Legendre duplication formula
\begin{equation*}
\forall p \in \nn \setminus \{0\}, \forall z \in \cc \setminus \left\{-\nn\right\}, \quad \prod \limits_{j=0}^{p-1} \Gamma \Big(\frac{z+j}{p} \Big) = (2 \pi)^{\frac{p-1}{2}} p^{\frac{1}{2}-z} \Gamma(z),
\end{equation*}
see e.g.~\cite[Chapter 3]{Artin}.
The Gamma function is strictly convex on $(0,+\infty)$, since differentiating under the integral sign provides that 
\begin{equation}
\forall x >0, \quad \Gamma''(x)= \int_{0}^{+\infty} (\ln t)^2 t^{x-1} e^{-t} dt >0.
\end{equation}
On the other hand, as $\Gamma(1)=\Gamma(2)=1$ thanks to (\ref{factoriel}), Rolle's theorem implies that there exists  $x_0$ in $\left]1,2\right[$ such that $\Gamma '(x_0)=0$. Since $\Gamma '$ is an increasing function on $(0,+\infty)$, the Gamma function is therefore increasing on $[2,+\infty)$. Related to the Gamma function is the Beta function
\begin{equation}
\forall x,y >0, \quad B(x,y)= \int_0^1 t^{x-1}(1-t)^{y-1} dt,
\end{equation}
satisfying the following identity 
\begin{equation}{\label{identitybeta}}
\forall x,y >0, \quad B(x,y)= \frac{\Gamma(x) \Gamma(y)}{\Gamma(x+y)}.
\end{equation}
Instrumental in the core of this work are the two following lemmas:

\medskip

\begin{lemma}\label{gamma4}
The Gamma function satisfies the following estimates:
\begin{equation*}
\forall p \in \nn \setminus \{0\},\exists C_p >0, \forall x \geq 1, \quad \Gamma(x)^{\frac{1}{p}} \leq C_p \big(p^{\frac{1}{p}}e^{\frac{1}{p}}\big)^x\Gamma\Big( \frac{x}{p} \Big).
\end{equation*}
\end{lemma}

\medskip

\begin{proof}
By using the Legendre duplication formula, (\ref{gamma2}) and the fact that the Gamma function is increasing on $[2,+\infty)$, we deduce that for all $p \in \nn \setminus \{0\}$, $x \geq 2p$,
\begin{multline*}
\Gamma(x) = (2 \pi)^{\frac{1-p}{2}} p^{x-\frac{1}{2}} \prod_{j=0}^{p-1} \Gamma \Big( \frac{x+j}{p} \Big) 
	\leq (2 \pi)^{\frac{1-p}{2}} p^{x-\frac{1}{2}} \Big(\Gamma \Big( \frac{x}{p} +1\Big)\Big)^p\\
	= (2 \pi)^{\frac{1-p}{2}} p^{x-\frac{1}{2}} \Big(\frac{x}{p} \Big)^p \Big(\Gamma\Big( \frac{x}{p}\Big)\Big)^p 
	\leq (2 \pi)^{\frac{1-p}{2}} p! p^{-\frac{1}{2}-p} \big(pe\big)^x \Big(\Gamma \Big( \frac{x}{p}\Big)\Big)^p,
\end{multline*}
since $x^p \leq e^x p!$. It proves the estimate when $x \geq 2p$.
We conclude by using the continuity of the function $x \longmapsto \frac{\Gamma(x)^{\frac{1}{p}}}{(p^{\frac{1}{p}}e^{\frac{1}{p}})^x \Gamma ( \frac{x}{p})}$ on $[1,+\infty)$.
\end{proof}

\medskip

\begin{lemma}{\label{gamma3}}
The Gamma function and the Beta function satisfy the following estimates:
\begin{align*}
(i) &\ \forall x >0, \forall y > 0, \quad x^y \leq \Gamma(y+1) e^x \\
(ii) &\ \forall r>0, \forall x,y \geq r, \quad \Gamma(x)\Gamma(y) \leq \frac{1}{2r}B(r,r) \Gamma(x+y+1) \\
(iii) &\ \forall r\geq1, \exists C_r>0, \forall x \geq 1, \quad \Gamma(x)^r \leq C_r\Gamma(rx)
\end{align*}
\end{lemma}

\medskip

\begin{proof}
It follows from $\eqref{gamma1}$ that for all $x, y>0$, 
\begin{equation*}
\Gamma(y)= \int_{0}^{+\infty} t^{y-1} e^{-t} dt  \geq \int_{0}^x t^{y-1} e^{-t} dt = x^{y} \int_{0}^1 t^{y-1} e^{-t x} dt 
\geq x^{y} \int_{0}^1 t^{y-1} e^{-x} dt = \frac{x^y e^{-x}}{y}.
\end{equation*}
Assertion $(i)$ directly follows from the previous estimate together with the functional identity $\eqref{gamma2}$.
On the other hand, since the Beta function is separately non-increasing with respect to the two variables, it follows from the functional identity $\eqref{gamma2}$ and (\ref{identitybeta}) that for all 
$r>0$, $x,y \geq r$, 
\begin{multline*}
\Gamma(x)\Gamma(y) = B(x,y)\Gamma(x+y)  \leq B(r,y)\Gamma(x+y)  \leq B(r,r) \Gamma(x+y) \\
=\frac{B(r,r)}{x+y} (x+y)\Gamma(x+y) 
\leq \frac{B(r,r)}{2r} \Gamma(x+y+1).
\end{multline*}
It proves that the estimate $(ii)$ holds. By using Stirling formula 
\begin{equation}\label{stirling1}
\Gamma(x) \sim_{x \to +\infty} \sqrt{\frac{2\pi}{x}} \Big( \frac{x}{e} \Big)^x, 
\end{equation}
see e.g.~\cite{Artin}, it follows that for all $r \geq 1$,
\begin{equation}
\frac{\Gamma(x)^r}{\Gamma(rx)} \sim_{x \to +\infty} \Big(\frac{2\pi}{x}\Big)^{\frac{r-1}{2}}r^{\frac{1}{2}-rx}= \mathcal{O}_r(1) \textrm{ when } x \to +\infty.
\end{equation}
Since the function $x \longmapsto \frac{\Gamma(x)^r}{\Gamma(rx)}$ is continuous on $[1,+\infty)$, there exists a positive constant $C_r >0$ such that the estimate $(iii)$ holds. 
\end{proof}

\subsection{Hermite functions and Bernstein type estimates}\label{weighted}
The standard Hermite functions $(\phi_{k})_{k\geq 0}$ are defined for $x \in \rr$,
 \begin{multline}\label{defi}
 \phi_{k}(x)=\frac{(-1)^k}{\sqrt{2^k k!\sqrt{\pi}}} e^{\frac{x^2}{2}}\frac{d^k}{dx^k}(e^{-x^2})
 =\frac{1}{\sqrt{2^k k!\sqrt{\pi}}} \Bigl(x-\frac{d}{dx}\Bigr)^k(e^{-\frac{x^2}{2}})=\frac{ a_{+}^k \phi_{0}}{\sqrt{k!}},
\end{multline}
where $a_{+}$ is the creation operator
$$a_{+}=\frac{1}{\sqrt{2}}\Big(x-\frac{d}{dx}\Big).$$
The Hermite functions satisfy the identity 
\begin{equation}\label{sd1}
\forall \xi \in \rr, \forall k \geq 0,  \quad \widehat{\phi_{k}}(\xi)=(-i)^k\sqrt{2\pi}\phi_k(\xi).
\end{equation}
The $L^2$-adjoint of the creation operator is the annihilation operator
$$a_-=a_+^*=\frac{1}{\sqrt{2}}\Big(x+\frac{d}{dx}\Big).$$
The following identities hold
\begin{equation}\label{eq2ui}
[a_-,a_+]=a_-a_+-a_+a_-=\text{Id}, \quad -\frac{d^2}{dx^2}+x^2=2a_+a_-+1,
\end{equation}
\begin{equation}\label{eq2}
\forall k \in \nn, \quad a_+ \phi_{k}=\sqrt{k+1} \phi_{k+1}, \qquad  \forall k \in \nn, \quad a_-\phi_{k}=\sqrt{k} \phi_{k-1} \ (=0\textrm{ si }k=0),
\end{equation}
\begin{equation}\label{eq2ui1}
\forall k \in \nn, \quad \Big(-\frac{d^2}{dx^2}+x^2\Big)\phi_{k}=(2k+1)\phi_{k}.
\end{equation}
The family $(\phi_{k})_{k\in \nn}$ is an orthonormal basis of $L^2(\R)$.
We set for $\alpha=(\alpha_{j})_{1\le j\le n}\in\N^n$, $x=(x_{j})_{1\le j\le n}\in \R^n,$
\begin{equation}\label{jk1}
\Phi_{\alpha}(x)=\prod_{j=1}^n\phi_{\alpha_j}(x_j).
\end{equation}
The family $(\Phi_{\alpha})_{\alpha \in \nn^n}$ is an orthonormal basis of $L^2(\R^n)$
composed of the eigenfunctions of the $n$-dimensional harmonic oscillator
\begin{equation}\label{6.harmo}
\mathcal{H}=-\Delta_x+|x|^2=\sum_{k\ge 0}(2k+n)\mathbb P_{k},\quad \text{Id}=\sum_{k \ge 0}\mathbb P_{k},
\end{equation}
where $\mathbb P_{k}$ is the orthogonal projection onto $\text{Span}_{\cc}
\{\Phi_{\alpha}\}_{\alpha\in \N^n,\val \alpha =k}$, with $\val \alpha=\alpha_{1}+\dots+\alpha_{n}$. 
Instrumental in the proof of Theorem~\ref{Spectral} are the following Bernstein type estimates:

\medskip

\begin{proposition}{\label{prop1}}
With $\mathcal E_{N}=\emph{\textrm{Span}}_{\cc}\{\Phi_{\alpha}\}_{\alpha \in \nn^n, \ |\alpha| \leq N}$, finite combinations of Hermite functions satisfy the following estimates: 
\begin{multline*}
\forall 0< \eps \leq 1, \exists K_\eps >1, \forall 0< \delta \leq 1, \exists \tilde{K}_{\eps,\delta} >1, \forall r >0, \forall \alpha, \beta \in \nn^n, \forall N \in \nn, \forall f \in \mathcal{E}_N,\\
\|x^\alpha \partial_x^\beta f \|_{L^2(\rr^n)} \leq \tilde{K}_{\eps,\delta} (\delta K_\eps)^{|\alpha|+|\beta|} \Gamma \Big(\frac{|\alpha|+|\beta|}{2-\eps}+2\Big) e^{\frac{N^{1-\frac{\eps}{2}}}{\delta^{2-\eps}}} \|f\|_{L^2(\rr^n)}, 
\end{multline*}
\begin{equation*}
 \|\left\langle x\right\rangle^r \partial_x^\beta f \|_{L^2(\rr^n)} \leq \tilde{K}_{\eps,\delta} K_\eps^{|\beta|+r} \delta^{|\beta|}\big( n+1 \big) ^{\frac{r}{2}} \Gamma \Big(\frac{r+|\beta|}{2-\eps}+3\Big) e^{\frac{N^{1-\frac{\eps}{2}}}{\delta^{2-\eps}}} \left\|f\right\|_{L^2(\rr^n)}.
\end{equation*}
\end{proposition}

\medskip

\begin{proof}
We notice that 
\begin{equation}\label{eq1}
x_j=\frac{1}{\sqrt{2}}(a_{j,+}+a_{j,-}), \quad \partial_{x_j}=\frac{1}{\sqrt{2}}(a_{j,-}-a_{j,+}),
\end{equation}
with 
\begin{equation}\label{ucla4}
a_{j,+}=\frac{1}{\sqrt{2}}(x_j-\partial_{x_j}), \quad a_{j,-}=\frac{1}{\sqrt{2}}(x_j+\partial_{x_j}).
\end{equation}
By denoting $(e_j)_{1 \leq j \leq n}$ the canonical basis of $\rr^n$, we obtain from (\ref{eq2}) and (\ref{eq1}) that for all $N \in \nn$ and $f \in \mathcal E_{N}$,
\begin{align*}
& \ \|a_{j,+}f\|_{L^2(\rr^n)}^2=\Big\|a_{j,+}
\Big(\sum_{|\alpha| \leq N}\langle f,\Phi_{\alpha}\rangle_{L^2}\Phi_{\alpha}\Big)\Big\|_{L^2(\rr^n)}^2\\
=& \ \Big\|\sum_{|\alpha| \leq N}\sqrt{\alpha_j+1}\langle f,\Phi_{\alpha}\rangle_{L^2}\Phi_{\alpha+e_j}\Big\|_{L^2(\rr^n)}^2
=\sum_{|\alpha| \leq N}(\alpha_j+1)|\langle f,\Phi_{\alpha}\rangle_{L^2}|^2\\
\leq & \ (N+1)\sum_{|\alpha| \leq N}|\langle f,\Phi_{\alpha}\rangle_{L^2}|^2=(N+1)\|f\|_{L^2(\rr^n)}^2
\end{align*}
and
\begin{align*}
& \ \|a_{j,-}f\|_{L^2(\rr^n)}^2=\Big\|a_{j,-}
\Big(\sum_{|\alpha| \leq N}\langle f,\Phi_{\alpha}\rangle_{L^2}\Phi_{\alpha}\Big)\Big\|_{L^2(\rr^n)}^2\\
=& \ \Big\|\sum_{|\alpha| \leq N}\sqrt{\alpha_j}\langle f,\Phi_{\alpha}\rangle_{L^2}\Phi_{\alpha-e_j}\Big\|_{L^2(\rr^n)}^2
=\sum_{|\alpha| \leq N}\alpha_j|\langle f,\Phi_{\alpha}\rangle_{L^2}|^2\\
\leq & \ N\sum_{|\alpha| \leq N}|\langle f,\Phi_{\alpha}\rangle_{L^2}|^2=N\|f\|_{L^2(\rr^n)}^2.
\end{align*}
It follows that for all $N \in \nn$ and $f \in \mathcal E_{N}$,
\begin{equation}\label{a1}
\|x_jf\|_{L^2(\rr^n)}\leq \frac{1}{\sqrt{2}}( \|a_{j,+}f\|_{L^2(\rr^n)}+ \|a_{j,-}f\|_{L^2(\rr^n)}) \leq \sqrt{2N+2}\|f\|_{L^2(\rr^n)}
\end{equation}
and
\begin{equation}\label{a2}
\|\partial_{x_j}f\|_{L^2(\rr^n)}\leq \frac{1}{\sqrt{2}}( \|a_{j,+}f\|_{L^2(\rr^n)}+ \|a_{j,-}f\|_{L^2(\rr^n)}) \leq \sqrt{2N+2}\|f\|_{L^2(\rr^n)}.
\end{equation}
We notice from (\ref{eq2}) and (\ref{eq1}) that
$$\forall N \in \nn, \forall f \in \mathcal E_{N}, \forall \alpha, \beta \in \nn^n, \quad x^{\alpha}\partial_x^{\beta}f \in  \mathcal E_{N+|\alpha|+|\beta|},$$
with $x^{\alpha}=x_1^{\alpha_1}...x_n^{\alpha_n}$ and $\partial_x^{\beta}=\partial_{x_1}^{\beta_1}...\partial_{x_n}^{\beta_n}$.
We deduce from (\ref{a1}) that for all $N \in \nn$, $f \in \mathcal E_{N}$, and $\alpha, \beta \in \nn^n$, with $\alpha_1 \geq 1$,
$$\|x^{\alpha}\partial_x^{\beta}f\|_{L^2(\rr^n)}=\|x_1(\underbrace{x^{\alpha-e_1}\partial_x^{\beta}f}_{\in \mathcal E_{N+|\alpha|+|\beta|-1}})\|_{L^2(\rr^n)}\leq
\sqrt{2}\sqrt{N+|\alpha|+|\beta|}\|x^{\alpha-e_1}\partial_x^{\beta}f\|_{L^2(\rr^n)}.$$
By iterating the previous estimates, we readily obtain from (\ref{a1}) and (\ref{a2}) that for all $N \in \nn$, $f \in \mathcal E_{N}$ and $\alpha, \beta \in \nn^n$, 
\begin{equation}\label{gh0}
\|x^{\alpha}\partial_x^{\beta}f\|_{L^2(\rr^n)}\leq 2^{\frac{|\alpha|+|\beta|}{2}}\sqrt{\frac{(N+|\alpha|+|\beta|)!}{N!}}\|f\|_{L^2(\rr^n)}.
\end{equation}
We recall the following basic estimates
\begin{equation}\label{rod1}
 \forall t \geq 0, \forall k \in \nn, \quad t^k \leq e^t k!, \quad \quad \forall t>0, \forall A>0,\quad  t^A \leq A^A e^{t-A}, 
\end{equation}
see e.g.~\cite{rodino} (formula (0.3.14)). 
Let $0<\delta \leq 1$ be a positive constant. When $N \leq |\alpha|+|\beta|$, with $|\alpha|+|\beta| \geq 1$, we deduce from (\ref{factoriel}) and (\ref{rod1}) that for all $p \in \nn \setminus \{0\}$, 
\begin{align}\label{gh1a}
& \ 2^{\frac{|\alpha|+|\beta|}{2}}\sqrt{\frac{(N+|\alpha|+|\beta|)!}{N!}} \leq 2^{\frac{|\alpha|+|\beta|}{2}}(N+|\alpha|+|\beta|)^{\frac{|\alpha|+|\beta|}{2}}\\ \notag
 \leq & \ 2^{|\alpha|+|\beta|}(|\alpha|+|\beta|)^{\frac{|\alpha|+|\beta|}{2}} \leq  (2\sqrt{e})^{|\alpha|+|\beta|}\sqrt{(|\alpha|+|\beta|)!}
  \\ = & \ \notag (2 \delta \sqrt{e})^{|\alpha|+|\beta|}\big(\Gamma(|\alpha|+|\beta|+1)\big)^{\frac{1}{2}} \Big(\frac{1}{\delta^p} \Big)^{\frac{|\alpha|+|\beta|}{p}} \\ \notag \leq   & \
(2 \delta \sqrt{e})^{|\alpha|+|\beta|}\big(\Gamma(|\alpha|+|\beta|+1)\big)^{\frac{1}{2}} \big((|\alpha|+|\beta|)! \big)^{\frac{1}{p}} e^{\frac{1}{p \delta^p}}
\\ \notag
=  & \  (2 \delta \sqrt{e})^{|\alpha|+|\beta|} e^{\frac{1}{p \delta^p}} \big(\Gamma(|\alpha|+|\beta|+1)\big)^{\frac{1}{2}}\big( \Gamma(|\alpha|+|\beta|+1)\big)^{\frac{1}{p}}.
\end{align}
The above estimate also holds when $|\alpha|+|\beta|=0$.
By using Lemma~\ref{gamma4} and Lemma~\ref{gamma3} (assertion $(ii)$), we deduce from $\eqref{gh1a}$ that for all $|\alpha|+|\beta| \geq N$, $0<\delta \leq 1$ and $p \in \nn \setminus \{0\}$,
\begin{align}{\label{gh1b}}
& \ 2^{\frac{|\alpha|+|\beta|}{2}}\sqrt{\frac{(N+|\alpha|+|\beta|)!}{N!}} \\  \notag
\leq & \  C_2 C_p e^{\frac{1}{p \delta^p}} \sqrt{2 e} e^{\frac{1}{p}} p^{\frac{1}{p}} \big(2\sqrt{2 }e e^{\frac{1}{p}} p^{\frac{1}{p}} \delta \big)^{|\alpha|+|\beta|} \Gamma \Big(\frac{|\alpha|+|\beta|+1}{2} \Big) \Gamma \Big(\frac{|\alpha|+|\beta|+1}{p} \Big)\\ \notag
\leq & \  \frac{p}{2}B\Big(\frac{1}{p},\frac{1}{p} \Big)C_2 C_p e^{\frac{1}{p \delta^p}} \sqrt{2 e} e^{\frac{1}{p}} p^{\frac{1}{p}} \big(2 \sqrt{2}e e^{\frac{1}{p}} p^{\frac{1}{p}} \delta\big)^{|\alpha|+|\beta|} \Gamma\Big( (|\alpha|+|\beta|)\Big(\frac{1}{2}+\frac{1}{p}\Big)+\Big(\frac{1}{2}+\frac{1}{p}+1\Big) \Big).
\end{align}
Let $0< \eps \leq 1$. We can choose the positive integer $p=p_{\eps}$ such that
\begin{equation}
p_\eps \geq 2, \quad \frac{1}{2}+\frac{1}{p_\eps} \leq \frac{1}{2-\eps}, \quad \frac{1}{2}+\frac{1}{p_\eps}+1 \leq 2.
\end{equation}
Since the Gamma function is convex on $(0,+\infty)$ and $\Gamma(1)=\Gamma(2)=1$, we have $\Gamma(x) \leq \Gamma(2)$ for all $1 \leq x \leq 2$. On the other hand, by using the fact that the Gamma function is increasing on $[2, +\infty)$, we deduce that $\Gamma(2) \leq \Gamma(y) \leq \Gamma(z)$ for all $2 \leq y \leq z$. It implies that 
\begin{equation}{\label{increasing2}}
\forall 1 \leq x \leq 2 \leq y \leq z, \quad \Gamma(x) \leq \Gamma(y) \leq \Gamma(z).
\end{equation}
It follows from (\ref{gh1b}) that $|\alpha|+|\beta| \geq N$, $0<\delta \leq 1$ and $0<\eps \leq 1$,
\begin{equation}{\label{gh1}}
2^{\frac{|\alpha|+|\beta|}{2}}\sqrt{\frac{(N+|\alpha|+|\beta|)!}{N!}} \leq \\
\tilde{K}_{\eps,\delta} \big(D_\eps \delta \big)^{|\alpha|+|\beta|} \Gamma \Big( \frac{|\alpha|+|\beta|}{2-\eps}+2 \Big),
\end{equation}
with
\begin{equation}
\tilde{K}_{\eps,\delta} =\frac{p_{\eps}}{2}B\Big(\frac{1}{p_{\eps}},\frac{1}{p_{\eps}}\Big) C_2 C_{p_\eps} e^{\frac{1}{p_\eps \delta^{p_\eps}}} \sqrt{2 e} e^{\frac{1}{p_\eps}} p_\eps^{\frac{1}{p_\eps}}>0 \text{ and  } D_\eps= 2 \sqrt{2}e e^{\frac{1}{p_\eps}} p_\eps^{\frac{1}{p_\eps}}>0.
\end{equation}
On the other hand, when $N \geq |\alpha|+|\beta|>0$, we deduce from (\ref{rod1}) and Lemma~\ref{gamma3} (assertion $(i)$) that for all $0<\delta \leq 1$, $0<\eps \leq 1$,
\begin{align}\label{gh2bis}
& \ 2^{\frac{|\alpha|+|\beta|}{2}}\sqrt{\frac{(N+|\alpha|+|\beta|)!}{N!}} \leq 2^{\frac{|\alpha|+|\beta|}{2}}(N+|\alpha|+|\beta|)^{\frac{|\alpha|+|\beta|}{2}}\\ \notag
\leq & \ (2\delta)^{|\alpha|+|\beta|}(\delta^{-1} \sqrt{N})^{|\alpha|+|\beta|}
=(2\delta)^{|\alpha|+|\beta|}(\delta^{\eps-2} N^{1-\frac{\eps}{2}})^{\frac{|\alpha|+|\beta|}{2-\eps}} \\ \notag
\leq & \ (2\delta)^{|\alpha|+|\beta|}\Big(\frac{|\alpha|+|\beta|}{2-\eps}\Big)^{\frac{|\alpha|+|\beta|}{2-\eps}}e^{\delta^{\eps -2}N^{1-\frac{\eps}{2}}-\frac{|\alpha|+|\beta|}{2-\eps}}\\ \notag
 \leq & \  (2\delta)^{|\alpha|+|\beta|}\Gamma \Big(\frac{|\alpha|+|\beta|}{2-\eps}+1 \Big )e^{\delta^{\eps -2}N^{1-\frac{\eps}{2}}}.
\end{align}
By using (\ref{increasing2}), we deduce from (\ref{gh2bis}) that for all $0<\delta \leq 1$, $0<\eps \leq 1$,
\begin{equation}\label{gh2}
2^{\frac{|\alpha|+|\beta|}{2}}\sqrt{\frac{(N+|\alpha|+|\beta|)!}{N!}} \leq (2\delta)^{|\alpha|+|\beta|}\Gamma \Big(\frac{|\alpha|+|\beta|}{2-\eps}+2 \Big )e^{\delta^{\eps -2}N^{1-\frac{\eps}{2}}},
\end{equation}
when $N \geq |\alpha|+|\beta|>0$. Let us also notice from (\ref{increasing2}) that 
\begin{equation}\label{asdf20}
2^{\frac{|\alpha|+|\beta|}{2}}\sqrt{\frac{(N+|\alpha|+|\beta|)!}{N!}} =1\leq (2\delta)^{|\alpha|+|\beta|}\Gamma \Big(\frac{|\alpha|+|\beta|}{2-\eps}+2 \Big )e^{\delta^{\eps -2}N^{1-\frac{\eps}{2}}},
\end{equation}
when $|\alpha|+|\beta|=0$, since $\Gamma(2)=1$.
It follows from (\ref{gh0}), (\ref{gh1}), (\ref{gh2}) and (\ref{asdf20}) that for all $0< \eps \leq 1$, there exists a positive constant $K_\eps > 1$ such that  
\begin{multline}\label{gh3}
\forall 0< \delta \leq 1, \exists \tilde{K}_{\eps,\delta} >1, \forall \alpha, \beta \in \nn^n, \forall N \in \nn, \forall f \in \mathcal{E}_N, \\ \|x^\alpha \partial_x^\beta f \|_{L^2(\rr^n)} \leq \tilde{K}_{\eps,\delta} (\delta K_\eps)^{|\alpha|+|\beta|} \Gamma \Big(\frac{|\alpha|+|\beta|}{2-\eps}+2\Big) e^{\frac{N^{1-\frac{\eps}{2}}}{\delta^{2-\eps}}} \left\|f\right\|_{L^2(\rr^n)}.
\end{multline}
By using Newton formula, we obtain that for all $k \in \nn$,
\begin{multline}\label{asdf21}
\|\left\langle x\right\rangle^k \partial_x^\beta f \|_{L^2(\rr^n)}^2 = \int_{\rr^n} \Big( 1 + \sum \limits_{i=1}^n {x_i^2} \Big)^k |\partial_x^\beta f(x) |^2 dx \\ = \int_{\rr^n} \sum_{\substack{\gamma \in \nn^{n+1}, \\ |\gamma|=k}} \frac{k!}{\gamma !} x^{2 \tilde{\gamma}} |\partial_x^\beta f(x) |^2 dx 
=\sum_{\substack{\gamma \in \nn^{n+1}, \\ |\gamma|=k}} \frac{k!}{\gamma !} \|x^{\tilde{\gamma}} \partial_x^\beta f \|_{L^2(\rr^n)}^2 ,
\end{multline}
where we denote $\tilde{\gamma}=(\gamma_1,...,\gamma_n) \in \nn^n$ if $\gamma=(\gamma_1,...\gamma_{n+1}) \in \nn^{n+1}$. It follows from (\ref{increasing2}), \eqref{gh3} and (\ref{asdf21}) that for all $0<\eps \leq 1$, $0<\delta \leq 1$, $\beta \in \nn^n$, $k \in \nn$, $N \in \nn$, $f \in \mathcal{E}_N$,
\begin{align}{\label{gh4}}
\|\left\langle x\right\rangle^k \partial_x^\beta f \|_{L^2(\rr^n)}^2 \leq & \  \sum_{\substack{\gamma \in \nn^{n+1}, \\ |\gamma|=k}} \frac{k!}{\gamma !} \tilde{K}_{\eps,\delta}^2 (\delta K_\eps)^{2|\tilde{\gamma}|+2|\beta|} \Big(\Gamma \Big(\frac{|\tilde{\gamma}|+|\beta|}{2-\eps}+2\Big)\Big)^2 e^{\frac{2N^{1-\frac{\eps}{2}}}{\delta^{2-\eps}}} \left\|f\right\|_{L^2(\rr^n)}^2 \\ \notag
\leq & \ \tilde{K}_{\eps,\delta}^2 \Big(\Gamma \Big(\frac{k+|\beta|}{2-\eps}+2\Big)\Big)^2 e^{\frac{2N^{1-\frac{\eps}{2}}}{\delta^{2-\eps}}} K_\eps^{2|\beta|+2k} \delta^{2|\beta|} \Big(\sum_{\substack{\gamma \in \nn^{n+1}, \\ |\gamma|=k}} \frac{k!}{\gamma !}\Big)  \left\|f\right\|_{L^2(\rr^n)}^2 \\ \notag
= & \ \tilde{K}_{\eps,\delta}^2 \Big(\Gamma \Big(\frac{k+|\beta|}{2-\eps}+2\Big)\Big)^2 e^{\frac{2N^{1-\frac{\eps}{2}}}{\delta^{2-\eps}}} K_\eps^{2|\beta|+2k}\delta^{2|\beta|}  (n+1)^k \left\|f\right\|_{L^2(\rr^n)}^2,
\end{align}
since 
\begin{equation}
\sum_{\substack{\gamma \in \nn^{n+1}, \\ |\gamma|=k}} \frac{k!}{\gamma !}=(n+1)^k,
\end{equation}
thanks to Newton formula. Let $r \in [0,+\infty) \setminus{\nn}$. We can write $r=\theta k+ (1-\theta) (k+1)>0$ with $k \in \nn$ and $\theta \in \left]0,1\right[$. By using H\"older inequality, it follows from $\eqref{gh4}$ that
\begin{multline}{\label{gh5}}
\|\left\langle x\right\rangle^r \partial_x^\beta f \|_{L^2(\rr^n)} \leq \|\langle x\rangle^k \partial_x^\beta f\|_{L^2(\rr^n)}^{\theta}\|\langle x\rangle^{k+1} \partial_x^\beta f\|_{L^2(\rr^n)}^{1-\theta} \\
\leq \tilde{K}_{\eps,\delta} \Big(\Gamma \Big(\frac{k+|\beta|}{2-\eps}+2\Big)\Big)^{\theta} \Big(\Gamma \Big(\frac{k+1+|\beta|}{2-\eps}+2\Big)\Big)^{1-\theta} e^{\frac{N^{1-\frac{\eps}{2}}}{\delta^{2-\eps}}} K_\eps^{|\beta|+ r} \delta^{|\beta|} (n+1)^{\frac{r}{2}} \|f\|_{L^2(\rr^n)}.
\end{multline}
By using that the Gamma function is increasing on $[2, +\infty)$ and that $k \leq r$, we deduce from $\eqref{gh5}$ that
\begin{align*}
& \ \|\left\langle x\right\rangle^r \partial_x^\beta f \|_{L^2(\rr^n)}\\ \notag
\leq & \  \tilde{K}_{\eps,\delta} \Big( \Gamma \Big(\frac{r+|\beta|}{2-\eps}+2\Big)\Big)^{\theta} \Big(\Gamma \Big(\frac{r+1+|\beta|}{2-\eps}+2\Big)\Big)^{1-\theta} e^{\frac{N^{1-\frac{\eps}{2}}}{\delta^{2-\eps}}} K_\eps^{|\beta|+ r}\delta^{|\beta|}  (n+1)^{\frac{r}{2}} \|f\|_{L^2(\rr^n)} \\ \notag
\leq & \ \tilde{K}_{\eps,\delta}  \Gamma \Big(\frac{r+|\beta|}{2-\eps}+3\Big) e^{\frac{N^{1-\frac{\eps}{2}}}{\delta^{2-\eps}}} K_\eps^{|\beta|+ r}\delta^{|\beta|}  (n+1)^{\frac{r}{2}} \|f\|_{L^2(\rr^n)},
\end{align*}
since $0< \frac{1}{2-\eps} \leq 1$, as $0< \eps \leq 1$. This ends the proof of Proposition~\ref{prop1}. 
\end{proof}

\subsection{Gelfand-Shilov regularity}\label{gelfand}

We refer the reader to the works~\cite{gelfand,rodino1,rodino,toft} and the references herein for extensive expositions of the Gelfand-Shilov regularity theory.
The Gelfand-Shilov spaces $S_{\nu}^{\mu}(\rr^n)$, with $\mu,\nu>0$, $\mu+\nu\geq 1$, are defined as the spaces of smooth functions $f \in C^{\infty}(\rr^n)$ satisfying the estimates
$$\exists A,C>0, \quad |\partial_x^{\alpha}f(x)| \leq C A^{|\alpha|}(\alpha !)^{\mu}e^{-\frac{1}{A}|x|^{1/\nu}}, \quad x \in \rr^n, \ \alpha \in \mathbb{N}^n,$$
or, equivalently
$$\exists A,C>0, \quad \sup_{x \in \rr^n}|x^{\beta}\partial_x^{\alpha}f(x)| \leq C A^{|\alpha|+|\beta|}(\alpha !)^{\mu}(\beta !)^{\nu}, \quad \alpha, \beta \in \mathbb{N}^n,$$
with $\alpha!=(\alpha_1!)...(\alpha_n!)$ if $\alpha=(\alpha_1,...,\alpha_n) \in \nn^n$.
These Gelfand-Shilov spaces  $S_{\nu}^{\mu}(\rr^n)$ may also be characterized as the spaces of Schwartz functions $f \in \mathscr{S}(\rr^n)$ satisfying the estimates
$$\exists C>0, \eps>0, \quad |f(x)| \leq C e^{-\eps|x|^{1/\nu}}, \quad x \in \rr^n, \qquad |\widehat{f}(\xi)| \leq C e^{-\eps|\xi|^{1/\mu}}, \quad \xi \in \rr^n.$$
In particular, we notice that Hermite functions belong to the symmetric Gelfand-Shilov space  $S_{1/2}^{1/2}(\rr^n)$. More generally, the symmetric Gelfand-Shilov spaces $S_{\mu}^{\mu}(\rr^n)$, with $\mu \geq 1/2$, can be nicely characterized through the decomposition into the Hermite basis $(\Phi_{\alpha})_{\alpha \in \mathbb{N}^n}$, see e.g. \cite[Proposition~1.2]{toft},
\begin{multline*}
f \in S_{\mu}^{\mu}(\rr^n) \Leftrightarrow f \in L^2(\rr^n), \ \exists t_0>0, \ \big\|\big(\langle f,\Phi_{\alpha}\rangle_{L^2}\exp({t_0|\alpha|^{\frac{1}{2\mu}})}\big)_{\alpha \in \mathbb{N}^n}\big\|_{l^2(\mathbb{N}^n)}<+\infty\\
\Leftrightarrow f \in L^2(\rr^n), \ \exists t_0>0, \ \|e^{t_0\mathcal{H}^{\frac{1}{2\mu}}}f\|_{L^2(\rr^n)}<+\infty,
\end{multline*}
where $\mathcal{H}=-\Delta_x+|x|^2$ stands for the harmonic oscillator.

\subsection{Slowly varying metrics}{\label{vsm}}
This section is devoted to recall basic facts about slowly varying metrics. We refer the reader to~\cite{Hormander} (Section 1.4) for the proofs of the following results.
Let $X$ be an open subset in a finite dimensional $\rr$-vector space $V$ and $\|\cdot\|_x$ a norm in $V$ depending on $x \in X$. The family of norms $(\|\cdot\|_x)_{x \in X}$ is said to define a slowly varying metric in $X$ if there exists a positive constant $C \geq 1$ such that for all $x \in X$ and for all $y \in V$ satisfying $\|y-x\|_x <1$, then $y \in X$ and 
\begin{equation}{\label{equiv}}
\forall v \in V, \quad \frac{1}{C} \|v \|_x \leq \|v\|_y \leq C \|v \|_x.
\end{equation}

\medskip

\begin{lemma}\label{slowmet}\cite[Example~1.4.8]{Hormander}.
Let $X$ be an open subset in a finite dimensional $\rr$-vector space $V$ and $d(x)$ a Lipschitz continuous function, positive in $X$ and zero in $V \setminus X$, satisfying 
\begin{equation*}
\forall x,y \in X, \quad |d(x) - d(y) | \leq \|x-y \|,
\end{equation*}
where $\|\cdot\|$ is a fixed norm in $V$. Then, the family of norms $(\|\cdot\|_x)_{x \in X}$ given by
\begin{equation*}
\|v\|_x= \frac{2 \|v\|}{d(x)}, \quad x \in X, v \in V,
\end{equation*}
defines a slowly varying metric in X.
\end{lemma}

\medskip

Let us consider the case when $X=V=\rr^n$ and $\|\cdot\|$ is the Euclidian norm. If $0 < \eps < 1$ and $0< R \leq \frac{1}{2(1-\eps)}$, then the gradient of the function $\rho_\eps(x)=R\left\langle x\right\rangle^{1-\eps}$ given by 
$$\forall x \in \rr^n, \quad \nabla \rho_\eps(x)=R(1-\eps) \frac{x}{\left\langle x\right\rangle^{1+\eps}},$$ 
satisfies $\| \nabla \rho_\eps\|_{L^{\infty}(\rr^n)} \leq \frac{1}{2}$. The mapping $\rho_{\eps}$ is then a $\frac{1}{2}$-Lipschitz positive function and Lemma~\ref{slowmet} shows that the family of norms  $\|\cdot\|_x= \frac{\|\cdot\|}{R \left\langle x\right\rangle^{1-\eps}}$ defines a slowly varying metric on $\rr^n$.

\medskip

\begin{theorem}{\label{slowmetric}}
\cite[Theorem~1.4.10]{Hormander}.
Let $X$ be an open subset in $V$ a $\rr$-vector space of finite dimension $n \geq 1$ and $(\|\cdot\|_x)_{x \in X}$ be a family of norms in $V$ defining a slowly varying metric. Then, there exists a sequence $(x_k)_{k \geq 0} \in X^{\nn}$ such that the balls
\begin{equation*}
B_k=\left\{x \in V:\ \|x-x_k \|_{x_k} <1 \right\} \subset X,
\end{equation*}
form a covering of $X$, 
$$X = \bigcup \limits_{k=0}^{+\infty} B_k,$$ 
such that the intersection of more than $N=\big(4 C^3+1 \big)^n$ two by two distinct balls $B_k$ is always empty, where $C \geq 1$ denotes the positive constant appearing in the slowness condition \emph{(\ref{equiv})}.
\end{theorem}

\subsection{Instrumental lemmas}\label{lemma_proof}

This section is devoted to the proofs of instrumental lemmas:

\medskip

\begin{lemma}\label{thick_comparison}
Let $\rho_1, \rho_2 : \rr^n \longrightarrow (0,+\infty)$ be two continuous positive functions satisfying 
$$\forall x \in \rr^n, \quad 0<\rho_1 (x) \leq \rho_2 (x).$$ 
If $\omega$ is a measurable subset of $\rr^n$ verifying 
\begin{equation}\label{thick_rho1}
\forall x \in \rr^n, \quad |\omega \cap B(x, \rho_1(x)) | \geq \gamma |B(x, \rho_1 (x))|,
\end{equation}
with $1- \frac{1}{6^n} < \gamma \leq 1$, 
where $B(y,r)$ denotes the Euclidean ball centered at $y \in \rr^n$ with radius $r>0$ and where $|A|$ denotes the Lebesgue measure of $A$,
then it satisfies 
\begin{equation}\label{thick_rho2}
\forall x \in \rr^n, \quad |\omega \cap B(x, \rho_2(x)) | \geq \tilde{\gamma} |B(x, \rho_2 (x))|,
\end{equation}
with $\tilde{\gamma}= 1-(1-\gamma)6^n >0$.
\end{lemma}

\medskip

\begin{proof}
Let $\omega$ be a measurable subset of $\rr^n$ satisfying (\ref{thick_rho1}) and $x \in \rr^n$.
We begin by recovering $\overline{B(x,\rho_2 (x))}$ by a finite number of balls $B\big(x_k, \frac{\rho_1 (x_k)}{3}\big)$ with $\rho_1 (x_k) \leq 3 \rho_2 (x_k)$. In order to do so, we first notice that $\overline{B(x,\rho_2 (x))}$ is a compact set and that 
\begin{equation}
\overline{B(x,\rho_2 (x))} \subset \bigcup \limits_{\substack{y \in \overline{B(x,\rho_2 (x))}, \\ \rho_1(y) \leq 3 \rho_2 (x)}} B\Big(y, \frac{\rho_1 (y)}{3}\Big).
\end{equation}
Indeed, if $y \in \overline{B(x,\rho_2 (x))}$ and $\rho_1 (y) > 3 \rho_2 ( x)$, then the continuous function defined for all $t \in [0,1]$ by $f(t)= \rho_1(ty+(1-t)x)$ satisfies $f(0)=\rho_1(x) \leq \rho_2(x)$ and $f(1)=\rho_1(y) >3 \rho_2(x)$. It follows that there exists $0 < t_0 < 1$ such that $\rho_1(z)=3 \rho_2(x)$ with $z= t_0 y+ (1-t_0)x \in B(x,\rho_2(x))$ and $y \in B(z, \frac{\rho_1(z)}{3})$, as 
$$\|z-x\|=t_0\|x-y\| <\rho_2(x), \qquad \|y-z\|=(1-t_0)\|x-y\| <\rho_2(x)=\frac{1}{3}\rho_1(z).$$
It follows that there exists a finite sequence $(x_{i_k})_{0 \leq k\leq N}$ of $\overline{B(x,\rho_2(x))}$ such that 
\begin{equation}\label{recov30}
\overline{B(x,\rho_2 (x))} \subset \bigcup \limits_{k=0}^N B\Big(x_{i_k}, \frac{\rho_1(x_{i_k})}{3}\Big)\quad  \text{ and }\quad  \forall 0 \leq k \leq N, \quad \rho_1 (x_{i_k}) \leq 3 \rho_2 (x).
\end{equation} 
We can now use the following covering lemma~\cite{rudin} (Lemma~7.3):

\medskip

\begin{lemma}[Vitali covering lemma] \label{rudin_recov}
Let $(y_i)_{0\leq i \leq N}$ be a finite sequence of $\rr^n$ and $(r_i)_{0 \leq i \leq N} \subset (0,+\infty)^{N+1}$. There exists a subset $S \subset \{0,...,N\}$ such that
\begin{itemize}
\item[$(i)$] The balls $(B(y_i, r_i))_{i \in S}$ are two by two disjoint 
\item[$(ii)$] $\bigcup \limits_{i=0}^N B(y_i, r_i) \subset \bigcup \limits_{i \in S} B(y_i, 3 r_i)$
\end{itemize}
\end{lemma}

\medskip

It follows from Lemma~\ref{rudin_recov} and (\ref{recov30}) that there exists a subset $S \subset \{0,...,N\}$ such that the balls $\big(B\big(x_{i_k}, \frac{\rho_1 (x_{i_k})}{3}\big)\big)_{k \in S}$ are two by two disjoint and satisfy
\begin{equation}\label{asdf10}
B(x,\rho_2 (x)) \subset \bigcup \limits_{k \in S} B(x_{i_k},\rho_1(x_{i_k})).
\end{equation}
We also notice that 
$$\bigsqcup \limits_{k \in S} B\Big(x_{i_k},\frac{\rho_1(x_{i_k})}{3}\Big) \subset B(x, 2 \rho_2( x)),$$ 
since, if $y \in B\big(x_{i_k},\frac{\rho_1(x_{i_k})}{3}\big)$ then
$$\|y-x\| \leq \|y-x_{i_k}\|+\|x_{i_k}-x\| <\frac{\rho_1(x_{i_k})}{3}+\rho_2(x) \leq 2 \rho_2(x).$$
It follows from (\ref{thick_rho1}) and (\ref{asdf10}) that
\begin{align*}
|\omega \cap B(x, \rho_2 (x))| = & |B(x,\rho_2(x))| - |(\rr^n \setminus \omega) \cap B(x,\rho_2(x))| \\ 
\geq & \ |B(x,\rho_2(x))| - \sum_{k \in S} |(\rr^n \setminus \omega) \cap B(x_{i_k},\rho_1(x_{i_k}))| \\
\geq & \  |B(x,\rho_2(x))| - \sum_{k \in S} (1-\gamma) | B(x_{i_k},\rho_1(x_{i_k}))| 
\end{align*}
and 
\begin{multline*}
|\omega \cap B(x, \rho_2 (x))| \geq |B(x,\rho_2(x))| - \sum_{k \in S} (1-\gamma)3^n \Big|B\Big(x_{i_k},\frac{\rho_1(x_{i_k})}{3}\Big)\Big| \\ = |B(x,\rho_2(x))| -(1-\gamma)3^n \Big|\bigsqcup \limits_{k \in S} B\Big(x_{i_k},\frac{\rho_1(x_{i_k})}{3}\Big)\Big| \geq |B(x,\rho_2(x))| -(1-\gamma)3^n | B(x,2\rho_2(x))|.
\end{multline*}
We deduce that 
\begin{equation*}
|\omega \cap B(x, \rho_2 (x))| \geq (1-(1-\gamma)6^n) |B(x,\rho_2(x))|.
\end{equation*}
This ends the proof of Lemma~\ref{thick_comparison}.
\end{proof}

\medskip

\begin{lemma}\label{equivv}
Let $\frac{1}{2}<s \leq 1$, $0<t_0 \leq 1$ and $A$ be a closed operator on $L^2(\rr^n)$ which is the  infinitesimal generator of a strongly continuous contraction semigroup $(e^{-tA})_{t \geq 0}$ on $L^2(\rr^n)$. If the estimates 
\begin{multline}\label{GS44}
\exists C > 1, \exists m_1>0, \exists m_2 \geq 0, \forall 0< t \leq t_0, \forall \alpha, \beta  \in \nn^n, \forall g \in L^2(\rr^n),\\
\| x^{\alpha} \partial_x^{\beta}( e^{-t A}g)\|_{L^{\infty}(\rr^n)} \leq  \frac{C^{1+|\alpha|+|\beta|}}{t^{m_1 (|\alpha|+|\beta|)+m_2}} (\alpha !)^{\frac{1}{2s}}(\beta!)^{\frac{1}{2s}}  \|g\|_{L^2(\rr^n)},
\end{multline}
hold, then the estimates 
\begin{multline}\label{GS55}
\exists \tilde{C} > 1, \exists \tilde{m}_1>0, \exists \tilde{m}_2 \geq 0, \forall 0< t \leq t_0, \forall \alpha, \beta  \in \nn^n, \forall g \in L^2(\rr^n),\\
\| x^{\alpha} \partial_x^{\beta}( e^{-t A}g)\|_{L^2(\rr^n)} \leq  \frac{\tilde{C}^{1+|\alpha|+|\beta|}}{t^{\tilde{m}_1 (|\alpha|+|\beta|)+\tilde{m}_2}} (\alpha !)^{\frac{1}{2s}}(\beta!)^{\frac{1}{2s}}  \|g\|_{L^2(\rr^n)},
\end{multline}
hold.
\end{lemma}

\medskip

\begin{proof}
We assume that the estimates (\ref{GS44}) hold. It follows that there exist some positive constants $(C_{\tilde{\alpha}}(n))_{\tilde{\alpha}\in \nn^n, |\tilde{\alpha}| \leq 2n}$ such that for all $0< t \leq t_0$, $\alpha, \beta  \in \nn^n$, $g \in L^2(\rr^n)$,
\begin{align}
& \ \| x^{\alpha} \partial_x^{\beta}( e^{-t A}g)\|_{L^2(\rr^n)} \\ \notag
\leq & \ \Big(\int_{\rr^n}\frac{dx}{(1+|x|^2)^{2n}}\Big)^{\frac{1}{2}}\|(1+|x|^2)^nx^{\alpha} \partial_x^{\beta}( e^{-t A}g)\|_{L^{\infty}(\rr^n)}\\ \notag
\leq & \ \sum_{\substack{\tilde{\alpha}\in \nn^n,\\ |\tilde{\alpha}| \leq 2n}}C_{\tilde{\alpha}}(n)\|x^{\alpha+\tilde{\alpha}} \partial_x^{\beta}( e^{-t A}g)\|_{L^{\infty}(\rr^n)} \\ \notag
\leq & \ \sum_{\substack{\tilde{\alpha}\in \nn^n,\\ |\tilde{\alpha}| \leq 2n}}C_{\tilde{\alpha}}(n)
\frac{C^{2n+1+|\alpha|+|\beta|}}{t^{m_1 (|\alpha|+|\beta|+2n)+m_2}} ((\alpha+\tilde{\alpha}) !)^{\frac{1}{2s}}(\beta!)^{\frac{1}{2s}}  \|g\|_{L^2(\rr^n)}.
\end{align}
By using from (\ref{fo1}) that 
\begin{multline*}
\forall \alpha, \tilde{\alpha}\in \nn^n,\ |\tilde{\alpha}| \leq 2n, \\
 (\alpha+\tilde{\alpha}) !=\prod_{j=1}^n(\alpha_j+\tilde{\alpha}_j) ! \leq \prod_{j=1}^n(\alpha_j+2n) !\leq \prod_{j=1}^n2^{\alpha_j+2n}(\alpha_j)!(2n) !\leq (4^n(2n)!)^n 2^{|\alpha|}\alpha!,$$
\end{multline*}
we obtain that the estimates (\ref{GS55}) hold with $\tilde{m}_1=m_1>0$ and $\tilde{m}_2=2n m_1+m_2>0$.
\end{proof}

\medskip

\begin{lemma}\label{dvpt}
For any $k \in \nn$, there exists a finite family of real numbers $(C_{l_1,l_2}^k)_{\substack{l_1, l_2 \in \nn, \\ 0 \leq l_1+ l_2 \leq k+1}}$ satisfying 
$$\prod_{j=0}^{k} \big((-1)^{j} \partial_{x} +x\big) = \sum_{\substack{l_1, l_2 \in \nn, \\ 0 \leq l_1+ l_2 \leq k+1}} C^k_{l_1,l_2} x^{l_1} \partial^{l_2}_{x}, \qquad x \in \rr,$$
and 
$$\forall l_1, l_2 \in \nn, \ 0 \leq l_1+ l_2 \leq k+1, \quad |C^k_{l_1, l_2}| \leq 3^k (k+1)^{\frac{k+1-l_1-l_2}{2}},$$
while using for short the following abusive notation for possibly non-commutative differential operators
$$\prod_{j=0}^kA_j(x,D_x):=A_0(x,D_x)...A_k(x,D_x).$$
\end{lemma}

\medskip

\begin{proof} 
We proceed by induction on $k \in \nn$. For $k=0$ or $k=1$, the result of Lemma~\ref{dvpt} readily holds. Let us assume that it holds true for $k \in \nn$. We observe that
\begin{align}\label{hp2}
& \ \prod_{j=0}^{k+1} \big((-1)^{j} \partial_x +x\big) =  \Big(\sum_{\substack{l_1, l_2 \in \nn, \\ 0 \leq l_1+ l_2 \leq k+1}} C^k_{l_1,l_2} x^{l_1} \partial^{l_2}_x\Big)\big((-1)^{k+1} \partial_x+x \big) \\ \notag
= & \ \sum_{\substack{l_1, l_2 \in \nn, \\ 0 \leq l_1+ l_2 \leq k+1}} (-1)^{k+1}C^k_{l_1,l_2} x^{l_1} \partial^{l_2+1}_x +C^k_{l_1,l_2} x^{l_1+1} \partial^{l_2}_x + l_2 C^k_{l_1,l_2} x^{l_1} \partial^{l_2-1}_x \\ \notag
=& \ \sum_{\substack{l_1 \geq 0,\ l_2 \geq 1,\\ 1 \leq l_1+ l_2 \leq k+2}} (-1)^{k+1}C^k_{l_1,l_2-1} x^{l_1} \partial^{l_2}_x + \sum_{\substack{l_1 \geq 1,\ l_2 \geq 0, \\ 1 \leq l_1 + l_2 \leq k+2}} C^k_{l_1-1,l_2} x^{l_1} \partial^{l_2}_x\\ \notag
&\  \qquad +\sum_{\substack{ \\ l_1, l_2 \in \nn, \\ 0 \leq l_1+l_2 \leq k}} (l_2+1) C^k_{l_1,l_2+1} x^{l_1} \partial^{l_2}_x.
\end{align}
By setting for all $l_1, l_2 \in \nn$ with $0 \leq l_1+l_2 \leq k+2$, 
\begin{multline}
C^{k+1}_{l_1, l_2} = (-1)^{k+1}C^k_{l_1,l_2-1} \un_{[1,+\infty)}(l_2)+C^k_{l_1-1,l_2} \un_{[1, +\infty)}(l_1) \\ + (l_2+1) C^k_{l_1,l_2+1} \un_{[0,k]}(l_1+l_2) ,
\end{multline}
we deduce from \eqref{hp2} that
\begin{equation}
\prod_{j=0}^{k+1} \big((-1)^{j} \partial_x +x\big) =\sum_{\substack{l_1,l_2 \in \nn, \\ 0 \leq l_1 + l_2 \leq k+2}} C^{k+1}_{l_1,l_2} x^{l_1} \partial^{l_2}_x.
\end{equation}
By using the induction property, we deduce that for all $l_1, l_2 \in \nn$ with $0 \leq l_1+l_2 \leq k+2$, 
\begin{align*}
& \ |C^{k+1}_{l_1,l_2}| \\
\leq & \ 3^k (k+1)^{\frac{k+1-l_1-(l_2-1)}{2}}+ 3^k (k+1)^{\frac{k+1-(l_1-1)-l_2}{2}} + 3^k (k+1)^{\frac{k+1-l_1-(l_2+1)}{2}+1} \un_{[0,k]}(l_1 +l_2)  \\
\leq & \ 3^k (k+2)^{\frac{k+2-l_1 -l_2}{2}} + 3^k (k+2)^{\frac{k+2-l_1-l_2}{2}} + 3^k (k+2)^{\frac{k+2-l_1-l_2}{2}} = 3^{k+1} (k+2)^{\frac{k+2-l_1 -l_2}{2}}.
\end{align*}
This ends the proof of Lemma~\ref{dvpt}.
\end{proof}

\end{document}